\documentclass[journal,onecolumn]{IEEEtran}
 
\usepackage[english]{babel}
\usepackage{graphicx}
\usepackage{multirow}
\usepackage{framed}
\usepackage[normalem]{ulem}
\usepackage{amsmath}
\usepackage{amsthm}
\usepackage{amssymb}
\usepackage{amsfonts}
\usepackage{enumerate}
\usepackage{xcolor}
\usepackage[utf8]{inputenc}

 \usepackage[bf]{subfigure}

\usepackage{enumitem}
\usepackage{fullpage}
\usepackage{float}
\usepackage{algorithm}
\usepackage{algpseudocode}
\usepackage{caption}

 \usepackage{lineno}

\theoremstyle{definition}

\newcommand{\sys}{HypOp}

\title{Distributed Constrained Combinatorial \emph{Op}timization leveraging \emph{Hyp}ergraph Neural Networks  \vspace*{-1pt}}

\usepackage{authblk}

\author[1,*]{Nasimeh Heydaribeni}
\author[1]{Xinrui Zhan}
\author[1]{Ruisi Zhang}
\author[2]{Tina Eliassi-Rad,}
\author[1]{Farinaz Koushanfar}
\affil[1]{Department of Electrical and Computer Engineering\\
  University of California, San Diego}
  
  \affil[2]{Khoury College of Computer Sciences, 
   Northeastern University, Boston, MA, USA\\

nheydaribeni, x5zhan, ruz032, fkoushanfar @ucsd.edu; t.eliassirad@northeastern.edu}

\affil[*]{Corresponding Author}

\date{}
\begin{document}
{\let\newpage\relax\maketitle}

\begin{abstract}
Scalable addressing of high dimensional constrained combinatorial optimization problems is a challenge that arises in several science and engineering disciplines. Recent work introduced novel application of graph neural networks for solving quadratic-cost combinatorial optimization problems. However, effective utilization of models such as graph neural networks to address general problems with higher order constraints is an unresolved challenge. This paper presents a framework, \sys{}, which  advances the state of the art for solving combinatorial optimization problems in several aspects: (i) it generalizes the prior results to higher order constrained problems with arbitrary cost functions by leveraging hypergraph neural networks; (ii) enables scalability to larger problems by introducing a new distributed and parallel training architecture; (iii) demonstrates generalizability across different problem formulations by transferring knowledge within the same hypergraph; (iv) substantially boosts the solution accuracy compared with the prior art by suggesting a fine-tuning step using simulated annealing; 
(v) shows a remarkable progress on numerous benchmark examples, including hypergraph MaxCut, satisfiability, and resource allocation problems, with notable  run time improvements using a combination of fine-tuning and distributed training techniques. We showcase the application of \sys{} in scientific discovery by solving a hypergraph MaxCut problem on NDC drug-substance hypergraph. Through extensive experimentation on various optimization problems, \sys{} demonstrates superiority over existing unsupervised learning-based solvers and generic optimization methods.
\end{abstract}

\addtocontents{toc}{\protect\iffalse}

\section{Introduction}
Combinatorial optimization is ubiquitous across science and industry. Scientists often need to make decisions about how to allocate resources, design experiments, schedule tasks, or select the most efficient pathways among numerous choices. Combinatorial optimization techniques can help in these situations by finding the optimal or near-optimal solutions, thus aiding in the decision-making process.
Furthermore, the integration of artificial intelligence (AI) into the field of scientific discovery is growing increasingly fluid, providing  means to enhance and accelerate research
\cite{wang2023scientific}.  
An approach to integrate AI into scientific discovery involves leveraging machine learning (ML) methods to expedite and improve the  combinatorial optimization techniques to solve extremely challenging  optimization tasks.  
Several  combinatorial optimization problems are proved to be NP-hard, rendering most existing solvers non-scalable.  Moreover, the continually expanding size of today's  datasets makes existing optimization methods inadequate for addressing constrained optimization problems on such vast scales.  
 To facilitate the development of scalable and rapid optimization algorithms, various learning-based approaches have been proposed in the literature \cite{schuetz2022combinatorial, cappart2021combinatorial}. 

Learning-based optimization methods can be classified into three main categories:  supervised learning, unsupervised learning, and reinforcement learning (RL). Supervised learning methods \cite{khalil2016learning, bai2019simgnn, gasse2019exact, nair2020solving, li2018combinatorial}  train a model to address  the given problems using a dataset of solved problem instances. However, these approaches exhibit limitations in terms of generalizability to problem types not present in the training dataset and tend to perform poorly on larger problem sizes. Unsupervised learning approaches \cite{schuetz2022combinatorial, karalias2020erdos, toenshoff2021graph} do not rely on datasets of solved instances. Instead, they train ML models using optimization objectives as their loss functions. Unsupervised methods offer several advantages over supervised ones, including enhanced generalizability and eliminating the need for datasets containing solved problem instances, which can be challenging to acquire. RL-based methods \cite{mirhoseini2021graph, yolcu2019learning, ma2019combinatorial, kool2018attention} hold promise for specific optimization problems, provided that lengthy training and fine-tuning times can be accommodated. 

Unsupervised learning based optimization methods can be  conceptualized as a fusion of gradient descent-based optimization techniques with learnable transformation functions. In particular, in an unsupervised learning-based optimization algorithm, the  optimization variables are computed through an ML model, and the objective function is optimized by applying gradient descent over the parameters of this model. In other words, instead of conducting gradient descent directly on the optimization variables, it is applied to the parameters of the  ML model responsible for generating them. The goal is to establish  more efficient optimization paths compared to conventional direct gradient descent. This approach can potentially facilitate better escapes from subpar local optima.
 Various transformation functions are utilized in the optimization literature  to render the optimization problems  tractable.  For example, there are numerous techniques to linearize or convexify an optimization problem by employing transformation functions that are often lossy \cite{asghari2022transformation}. 
We believe that adopting ML-based transformation functions for optimization offers a substantial capability to train  customized transformation functions that best suit the specific optimization task.

 For combinatorial optimization problems over graphs, graph neural networks are  frequently employed as the learnable transformation functions \cite{schuetz2022combinatorial, karalias2020erdos, toenshoff2021graph}. However, when dealing with complex systems with higher order interactions and constraints, graphs fall short of modeling such relations. When intricate relationships among multiple entities extend beyond basic pairwise connections,  
  \emph{hypergraphs} emerge as invaluable tools for representing a diverse array of scientific phenomena. 
 They have been utilized in various areas including biology and bioinformatics \cite{feng2021hypergraph2, murgas2022hypergraph}, social network analysis \cite{zhu2018social}, chemical and material science \cite{xia2022novel}, image and data processing \cite{wen2012hyperspectral}. Hypergraph neural networks (HyperGNN) \cite{feng2019hypergraph} have also been commonly used for certain learning tasks such as image and visual object classification \cite{feng2019hypergraph},  and material removal rate prediction \cite{xia2022novel}. It is anticipated that HyperGNNs may serve as valuable tools for addressing  combinatorial optimization problems with higher-order constraints.

 While unsupervised learning approaches for solving combinatorial optimization problems offer numerous advantages, they may face challenges, including susceptibility to getting trapped in suboptimal local solutions and scalability issues.
 In a recent work \cite{angelini2022modern}, the authors argue that for certain well-known combinatorial optimization problems, unsupervised learning optimization methods may exhibit inferior performance compared to straightforward heuristics.
 However, it is crucial to recognize that these unsupervised methods possess a notable advantage in their generic applicability. Generic solvers like gradient-based techniques (such as SGD and ADAM), as well as simulated annealing (SA), may not be able to compete with customized heuristics that are meticulously crafted for specific problems. Nonetheless, they serve as invaluable tools for addressing less-known problems lacking effective heuristics. Consequently, the strategic utilization of unsupervised learning-based optimization methods can enhance and extend the capabilities of existing generic solvers, leading to the development of  efficient tools for addressing a diverse range of optimization problems.

In this study, we build upon the unsupervised learning-based optimization method for quadratic-cost combinatorial optimization  problems on graphs introduced in \cite{schuetz2022combinatorial}, and present \sys{}, a new scalable solver for a wide range of constrained combinatorial optimization problems with arbitrary cost functions. Our approach is applicable to problems with higher-order constraints by adopting hypergraph modeling for such problems (Fig.\,1(a)); \sys{} subsequently utilizes hypergraph neural networks in the training process,  a generalization of the graph neural networks employed in \cite{schuetz2022combinatorial}. 
\sys{} further  combines unsupervised HyperGNNs with another generic optimization method,  simulated annealing \cite{kirkpatrick1983optimization} to boost its  performance. 
Incorporating SA with HyperGNN can help with mitigation of the potential subpar local optima that may arise from HyperGNN training. 

To  establish a scalable solver and expedite the optimization process, \sys{} proposes two algorithms for parallel and distributed training of HyperGNN. First, it develops a distributed training algorithm in which, the hypergraph is distributed across a number of servers and each server only has a local view of the hypergraph. 
We develop a collaborative distributed algorithm to train the HyperGNN and solve the optimization problem (See Fig.\,1(b)). Second, \sys{} proposes a parallel HyperGNN training approach where the costs associated to constraints are computed in a scalable manner.  We further exhibit the transferability of our models, highlighting their efficacy in solving different optimization problems on the same hypergraph through transfer learning. This not only shows the generalizability of \sys{} but also considerably accelerates the optimization process.
 \sys{} is tested by comprehensive experiments, thoughtfully designed to provide insights into unsupervised learning-based optimization algorithms and their effectiveness across diverse problem types. We validate the scalability of \sys{} by testing it on  huge graphs, showcasing how parallel and distributed training can yield high-quality solutions even on such massive graphs. In summary, our contributions are as follows. 
\begin{figure}
    \centering
   \subfigure[Hypergraph Modeling of the Constraint Combinatorial Optimization Problem in \sys{}.] {\label{fig:cons_hyp}\includegraphics[width=5.9cm]{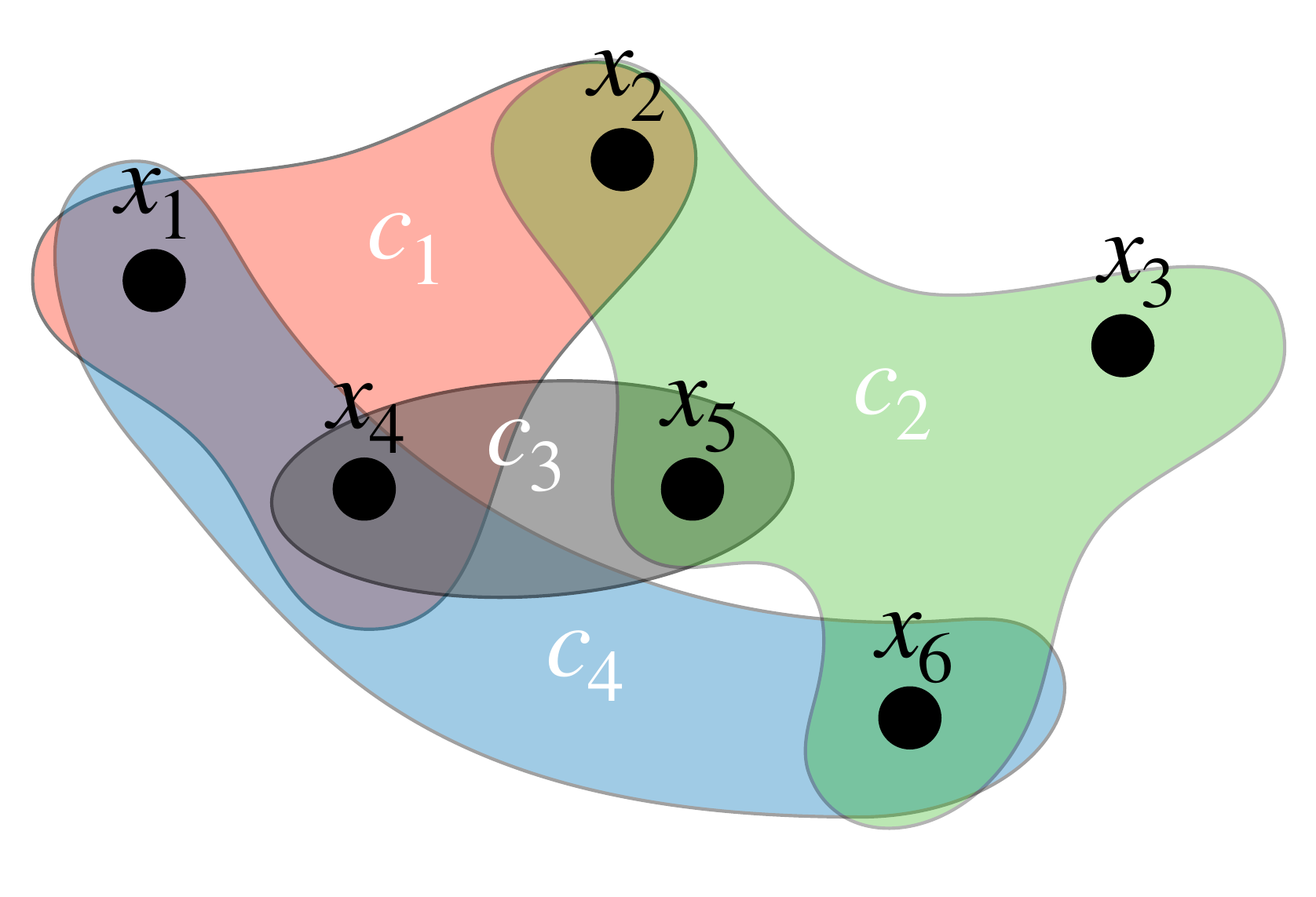}}
     \subfigure[Distributed HyperGNN training in \sys{}.]{\label{fig:dist:train}\includegraphics[width=7cm]{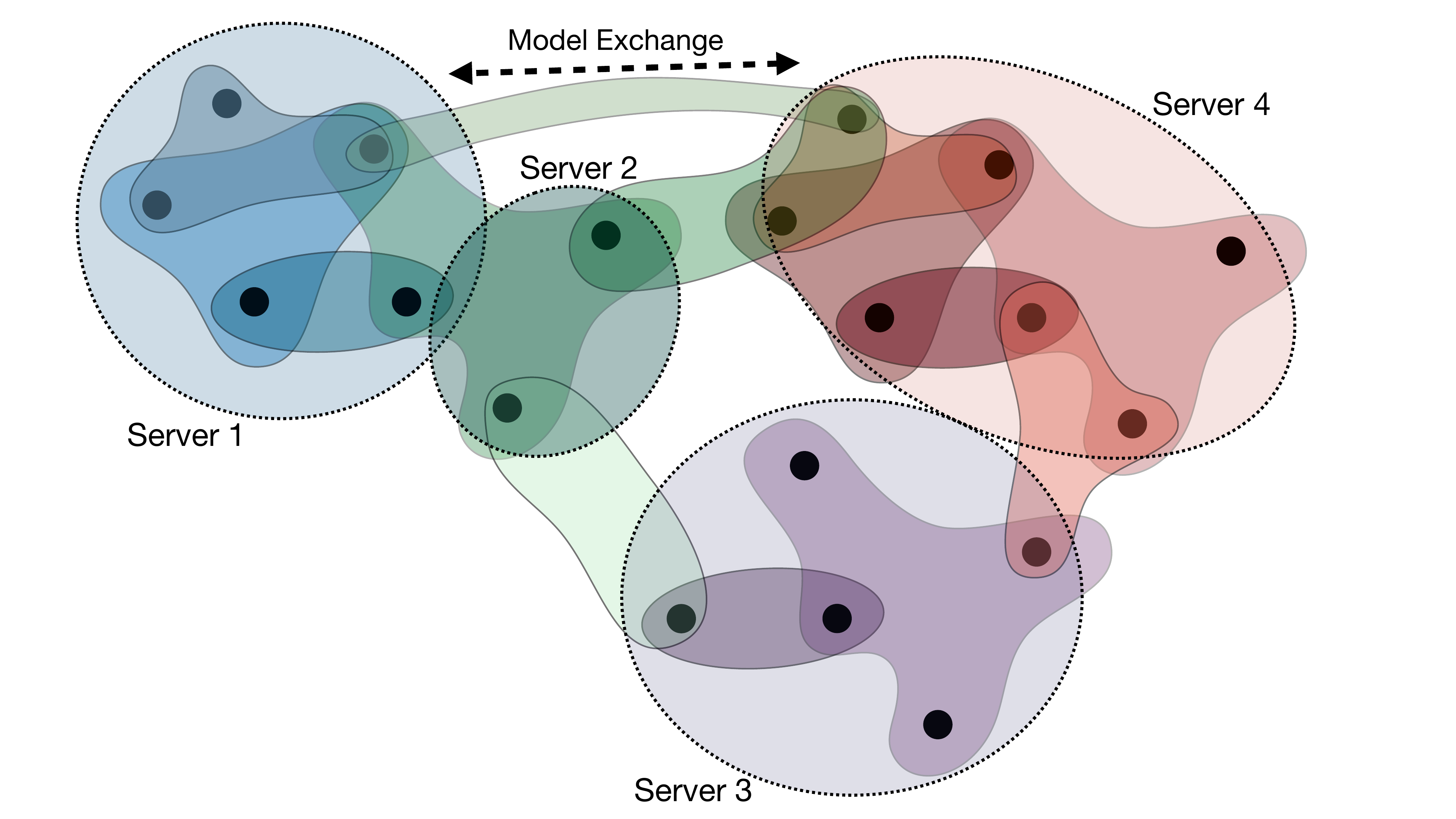}}
    \caption{Hypergraph Modeling and Distributed Training of HyperGNN in \sys{}}
\end{figure}
 \begin{itemize}
 \item Presenting \sys{}, a scalable unsupervised learning-based optimization method for addressing a wide spectrum of constrained combinatorial optimization problems with arbitrary cost functions. Notably, we pioneer the application of hypergraph neural networks within the realm of learning-based optimization for general  combinatorial optimization problems with higher order constraints. 
\item  Enabling scalability to much larger problems by introducing a new distributed and parallel architecture for hypergraph neural network training.
\item Demonstrating generalizability to other problem formulations, by knowledge transfer from the learned experience of addressing one set of  cost/constraints to another set for the same hypergraph.
\item Substantially boosting the solution accuracy and improving run time using a combination of fine-tuning (simulated annealing) and distributed training techniques.   
\item Demonstrating superiority of \sys{} over existing unsupervised learning-based solvers and generic optimization methods through extensive experimentation on a variety of combinatorial optimization problems. We address a novel set of scientific problems including hypergraph MaxCut problem, satisfiability problems (3SAT), and resource allocation. 
\item Showcasing  the application of \sys{} in scientific discovery by solving a hypergraph MaxCut problem on the NDC drug-substance hypergraph.
    \end{itemize}

\subsection{Paper Structure}
The remainder of this paper is structured as follows.   We provide the problem statement in Section \ref{sec:problem}. The \sys{} method is presented in Section \ref{sec:hypop}. We describe two algorithms for distributed and scalable training of \sys{} in Section \ref{sec:multigpu}. Our experimental results are provided in Section \ref{sec:exp}. We showcase the possibility of transfer learning in \sys{} in Section \ref{sec:transfer}. The applications of \sys{} in scientific discovery is discussed in Section \ref{sec:science}.
We conclude by a discussion in Section \ref{sec:discussion}. Supplementary information including, related work,  preliminaries, extra details on the experimental setup and results, and the limitations are provided in the Supplementary Information section. 
    
\section{Problem Statement}
\label{sec:problem}
We consider a constrained combinatorial optimization problem over $x_i,\ i \in \mathcal{N}$, where $\mathcal{N}=\{1, \cdots, N\}$ is the set of optimization variable indices. The variables $x_i$ are integers belonging to the set $\{d_0,\cdots, d_v\}$, where $d_0, \cdots, d_v$ are some given integers and $d_i<d_j$ for $i<j$. There are $K$ constraint functions, $c_k(x_{\mathcal{N}_k}), \ k \in \mathcal{K}$, where $\mathcal{K}=\{1, \cdots, K\}$ is the set of constraint indices, $x_{\mathcal{N}_k}=\{x_i, i \in \mathcal{N}_k\}$, and we have $\mathcal{N}_k\subset \mathcal{N}$. We consider the following optimization problem with an arbitrary cost function $f(x)$. 
\begin{subequations}
\begin{align}
    \min_{x_i, i \in \mathcal{N}} \quad &f(x)\\
    s.t. \quad &c_k(x_{\mathcal{N}_k})\leq 0, \ \text{for}\   k\in \mathcal{K}\\
 &x_i \in \{d_0,\cdots, d_v\}, \ \text{for}\ i\in \mathcal{N}
\end{align}
\label{eq:problem}
\end{subequations}
 Solving the above optimization problem  can be challenging, potentially falling into the category of NP-hard problems. This complexity arises from the presence of discrete optimization variables, as well as the potentially non-convex and non-linear  objective and constraint functions. Consequently, efficiently solving  this problem  is a complex task. Additionally, the discreteness of optimization variables implies that the objective function is not continuous, making gradient descent-based approaches impractical. To address this, a common strategy involves relaxing the state space into continuous spaces, solving the continuous problem, and then mapping the optimized continuous values back to discrete ones. While this approach may not guarantee globally optimal solutions due to the relaxation and potential convergence to suboptimal local optima, it remains a prevalent practice due to the availability and efficiency of the gradient descent-based solvers.  Such solvers, such as ADAM \cite{kingma2014adam}, are widely developed and used as the main optimization method for ML tasks. The relaxed version of the optimization problem \eqref{eq:problem} is given below. 
\begin{subequations}
\begin{align}
    \min_{p_i, i \in \mathcal{N}} \quad &f(p)\\
    s.t. \quad &c_k(p_{\mathcal{N}_k})\leq 0, \ \text{for}\   k\in \mathcal{K}\\
 &p_i \in [d_0, d_v], \ \text{for}\ i\in \mathcal{N},
\end{align}
\label{eq:problem:cont}
\end{subequations}
where  $p_i$ represents the continuous (relaxed) form of $x_i$, and it is within the continuous interval  $[d_0,d_v]$.

\section{\sys{} Method}
\label{sec:hypop}
In this section, we describe our method, called \sys{}, for solving problem \eqref{eq:problem}. An overview of our algorithm is provided in Fig.\,2, with each component explained subsequently. The main steps of the algorithm are: (1) hypergraph modeling of the problem, (2) solving the relaxed version of the problem (problem \eqref{eq:problem:cont}) using hypergraph neural networks in an unsupervised learning manner, and (3) mapping the continuous outcomes to discrete values using a post-processing mapping algorithm. 
\begin{figure}[t!]
\hspace{-0.45cm}
    \includegraphics[width=17cm]{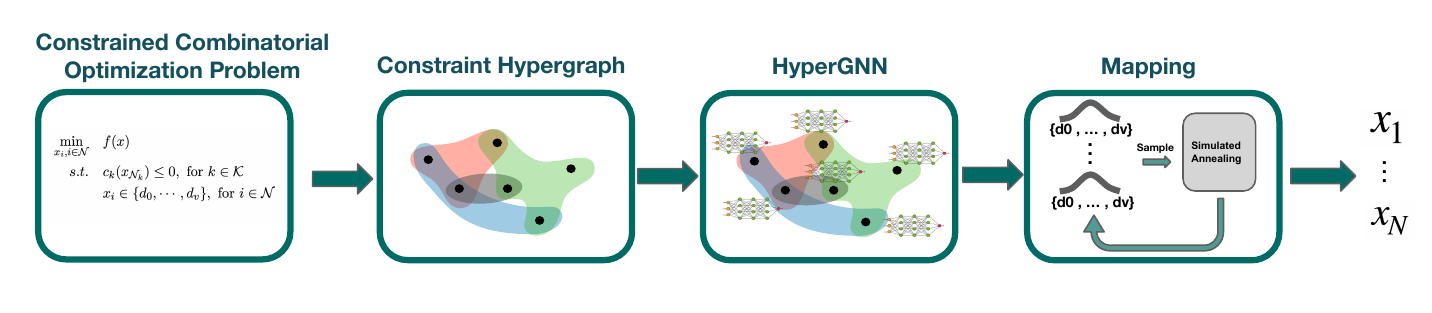}
    \caption{The overview of \sys{} for solving problem \eqref{eq:problem}.}
    \label{fig:overview}
\end{figure}
In the following, we describe each step of \sys{} in more details.

\subsection{Hypergraph Modeling}
\label{sec:cons_hyp}
To leverage hypergraph neural networks to solve problem \eqref{eq:problem} (or \eqref{eq:problem:cont}), we model the constraints of this problem through a hypergraph as follows. Consider a hypergraph with nodes $x_i, \ i \in \mathcal{N}$. For each constraint $c_k(x_{\mathcal{N}k})$, we form a hyperedge comprising nodes $x_{\mathcal{N}_k}$. This hypergraph is referred to as the constraint hypergraph of the problem. 
 For instance, in a problem with 6 variables and 4 constraints, the constraint hypergraph is illustrated in Fig.\,1(a). It is worth noting that for problems with an inherent graph or hypergraph structure, such as graph MaxCut, the constraint hypergraph coincides with the underlying graph or hypergraph of the problem.

\subsection{HyperGNN as the Learnable Transformation Function}
\label{sec:transform }
One can solve  optimization problem \eqref{eq:problem:cont} directly using gradient descent-based methods and then map the continuous outcomes to the discrete values. Penalty-based methods are a common practice in these approaches to ensure the constraints are satisfied in the final outcome. However, as will be  demonstrated in this paper, one can achieve better solutions  by solving an alternative optimization problem. In this alternative problem, the optimization variables $p_i$ are the outcomes of a transformation function based on HyperGNNs. Such transformation functions are utilized with the goal of capturing the intricate patterns and interdependence of the variables  through the problem constraints and therefore, obtaining better solutions.

In optimization literature, various transformation functions are utilized to render optimization problems tractable. For instance, there are numerous techniques for linearizing optimization problems through the use of transformation functions \cite{asghari2022transformation}. However, in much of the existing literature, the transformation function is a fixed function. In \sys{}, we introduce a parameterized transformation function based on HyperGNNs. 
 In particular, we consider $p=\mathcal{G}(\sigma;W)$ (or $p_i=\mathcal{G}_i(\sigma;W)$ for $i \in \mathcal{N}$), where $p_i$ is the (continuous) assignment to node $i$, $W$ is the parameter vector of the HyperGNN-based transformation function, $\mathcal{G}$, and $\sigma=[\sigma_1, \cdots, \sigma_N]$ is the input embedding ($\sigma_i$ is node $i$'s input embedding). 
By adopting a parameterized transformation function, we optimize our objective by simultaneously learning a good transformation function ($W$) and  optimizing the input embedding ($\sigma$). 
This alternative optimization problem is described below.
\begin{subequations}
\begin{align}
    \min_{\sigma, W} \quad &f(\mathcal{G}(\sigma;W))\\
    s.t. \quad &c_k(\mathcal{G}_{\mathcal{N}_k}(\sigma;W))\leq 0, \ \text{for}\   k\in \mathcal{K}\\
 &\mathcal{G}_i(\sigma;W) \in [d_0, d_v], \ \text{for}\ i\in \mathcal{N}
\end{align}
\label{eq:problem:trans}
\end{subequations}
We define the HyperGNN $\mathcal{G}(\sigma;W)$ on the constraint hypergraph with the following augmented (penalized)  loss function.
\begin{align}
    \hat{f}(p)=f(p)+\sum_{k \in {\mathcal{K}}}\lambda_k(c_k(p_{\mathcal{N}_k}))^+
    \label{eq:\sys{}_loss}
\end{align}
where $(a)^+=\max(0,a)$, and $\lambda_k$ is the weight of constraint $k$. 
The HyperGNN model used in \sys{} is based on hypergraph convolutional networks (HyperGCN)  introduced in \cite{feng2019hypergraph}. 
In particular, we have the following layer-wise operation in our HyperGNN model. 
 \begin{align}
    H^{(l+1)}=\sigma(D_v^{-\frac{1}{2}}\tilde{M}D_v^{-\frac{1}{2}}H^{(l)}W^{(l)})
\end{align}
where we define 
\begin{align}
\tilde{M}=A\tilde{D}_e^{-1}A^{T}- diag(A\tilde{D}_e^{-1}A^{T})
\end{align}
 and $H^{(l)}\in \mathbb{R}^{N\times f_l}$ is the matrix of node features at $l_{th}$ layer, $f_l$ is the size of the $l_{th}$ layer node features, $A$ is the hypergraph incidence matrix defined as $(A)_{ij}=1$ if node $i$ belongs to hyperedge $j$, and $(A)_{ij}=0$ otherwise. $D_v$ and $D_e$ are the diagonal degree matrices of nodes and hyperedges, respectively, and $\tilde{D}_e=D_e-I$.
  $W^{(l)}\in \mathbb{R}^{f_l\times f_{l+1}}$ is the $l_{th}$ layer  trainable weight matrix.

\subsection{Mapping with Simulated Annealing}
In order to transform the continuous outputs of the HyperGNN to integer node assignments, we use the continuous outputs generated through HyperGNN to derive a probability distribution over the discrete values  $\{d_0, \cdots, d_v\}$. This probability distribution is constructed such that the probability of a variable being equal to a given value is inversely proportional to the distance of the variable's continuous output with that value. For example, for the binary case of 0 and 1, the output will be the probability of 1. 
The mapping is then proceeded as follows. We take a sample, $x_1,\cdots, x_N$ from the output distribution corresponding to $p_1, \cdots, p_N$. We then initialize a simulated annealing algorithm with that  sample  and minimize  $\hat{f}(x)$. This process is repeated for a given number of times. The best outcome is chosen as the final solution.

\section{Distributed and Parallel \sys{} Training}
\label{sec:multigpu}
In this section, we introduce two scalable multi-GPU HyperGNN training algorithms for \sys{}: \textbf{Parallel Multi-GPU Training} that shuffles constraints among multiple servers (GPUs) and trains the HyperGNN in a parallel way; and \textbf{Distributed Training} where each server (GPU) holds part of the hypergraph (local view of the hypergraph) and trains the HyperGNN in a distributed  manner in collaboration with other servers. For a comparison between our approach  and the existing state of the art, see Supplementary Information Section I.

\subsection{Distributed Training}
In distributed training, we assume that the hypergraph is either inherently distributed across a number of servers, and these servers collaborate to solve the optimization problem, or the hypergraph is partitioned between the GPUs (servers) beforehand to expedite the optimization process and handle large-scale problems effectively.  
In this setup, each GPU  maintains a local view of the hypergraph (has access to a subgraph of the original hypergraph). The training of \sys{} takes place collaboratively among these GPUs.
Fig.\,1(b) shows a schematic example of distributed HyperGNN training in \sys{}.  

The outline of the distributed multi-GPU training is summarized in Algorithm\,\ref{alg:decentrailized_gpu}. Given a hypergraph $G = (V, E)$, and $S$ GPUs to train the HyperGNN model $M$, we first partition the hypergraph $G$ into $S$ subgraphs. The subgraph $G^{s} = (V^{s}, E^{s,\text{inner}}, E^{s, \text{outer}})$  consists of three components: $V^{s}$ represents the nodes in $G^{s}$, $E^{s,\text{inner}}$ includes inner hyperedges connecting nodes within $G^{s}$, and $E^{s, \text{outer}}$ consists of outer hyperedges linking nodes in $G^{s}$ to other subgraphs. 
During the forward pass, each GPU computes local node embeddings using $E^{s,\text{inner}}$, and during the backward pass, it calculates the loss with both local node embeddings and node embeddings connected by $E^{s,\text{outer}}$ to update model weights $M^s$. Subsequently, the gradients of each $M^s$ are aggregated and disseminated to all GPUs.

\begin{algorithm}[!ht]
\caption{Distributed Multi-GPU HyperGNN Training}
\label{alg:decentrailized_gpu}
\begin{algorithmic}[1]
\Require Hypergraph $G = (V, E)$, HyperGNN Model $M$, Number of GPUs $S$, Epoch Number $EP$
\Ensure Trained HyperGNN Model $M$

\State Partition $G$ into $S$ subgraphs: $G^{1}, G^{2}, \ldots, G^{S}$
\State Initialize $M$ on each GPU
\State Distribute $G^{s}$ to GPU$^s$, for $s=1$ to $S$

\For{$ep$ in $EP$}
    \For{$s \leftarrow 1$ to $S$}
        \State Extract inner hyperedges $E^{s,\text{inner}}$ and outer hyperedges $E^{s,\text{outer}}$ from $G^{s}$
        
        \State \textbf{Forward Pass:} Compute embeddings using $E^{s,\text{inner}}$ on GPU$^s$
        \State \textbf{Backward Pass:} Calculate loss using embeddings, $E^{s,\text{inner}}$, and $E^{s,\text{outer}}$
        
        \State Compute gradients based on loss
    \EndFor
    
    \State Aggregate gradients from all GPUs
    \State Update $M$ based on aggregated gradients
    \State Broadcast updated $M$ to all GPUs
\EndFor

\State \Return $M$

\end{algorithmic}
\end{algorithm}

\subsection{Parallel Multi-GPU Training}
HyperGNN-based optimization tool used in \sys{} allows us to take advantage of GNN and HyperGNN parallel training techniques to accelerate the optimization and enable optimization of large-scale problems. 
The outline of the parallel multi-GPU training is summarized in Supplementary Information Algorithm\,2. Given an input Hypergraph $G = (V, E)$, and $S$ GPUs,  the HyperGNN model $M$ is initialized on each GPU. In each epoch, the Hypergraph $G$'s hyperedge information $E$ are randomly partitioned into $S$ parts. Then each GPU$^s$ receives its corresponding hyperedges $E^s$ and trains $M^s$ on them. The gradients  are then aggregated and  disseminated to all GPUs.

\section{Experiments}
\label{sec:exp}

We have tested \sys{} on multiple types of combinatorial optimization problems, including  Hypergraph and Graph MaxCut, Graph Maximum Independent Set (MIS), Satisfiability (SAT), and Resource Allocation problems. We note that we have tested our method on the known and highly studied problems of SAT,  graph MaxCut, and MIS to show that even though our solver might not be able to compete with heuristics that are specifically designed for those problems, it can still be used as a descent solver for such well known problems.  Some of our results are provided in the Supplementary Information Section VII.

\subsection{Baselines}
Since \sys{} is proposed as a generic solver for combinatorial optimization problems, we compare \sys{} with other generic solvers such as simulated annealing (SA) \cite{kirkpatrick1983optimization} or gradient descent-based methods (e.g., ADAM) \cite{kingma2014adam}. Specifically, we use the following optimization methods as our baselines:

\textbf{- SA:} The experiments with SA solver are conducted by setting the HyperGNN training epochs in \sys{} to 0 and randomly initializing the SA step. 

\textbf{- ADAM:} The experiments with ADAM solver are conducted by removing the HyperGNN part from \sys{} and directly optimizing w.r.t. the continuous variables $p_i$ using ADAM. The exact SA fine-tuning step in \sys{} is then applied to the optimized results obtained using ADAM to generate discrete variable assignments.

\textbf{- Bipartite GNN:} In order to justify the need for HyperGNNs to solve problems with higher order constraints, we compare \sys{} with a GNN-based solver similar to PI-GNN \cite{schuetz2022combinatorial}. This GNN-based solver is developed by transforming the constraint hypergraph of the problem into a constraint Bipartite graph, where each hyperedge is represented by a factor node, which is connected to the nodes within that hyperedge. Then, \sys{} solves the problem on the bipartite graph. Note that this is similar to PI-GNN being applied to a bipartite graph, with the addition of the fine-tuning step that is provided in \sys{}.

\textbf{- PI-GNN:} For problems over graphs, such as graph MIS and MaxCut, we compare with PI-GNN \cite{schuetz2022combinatorial}.

\textbf{- Run-CSP:} For the MaxCut problem over graphs, we also consider Run-CSP \cite{toenshoff2021graph} as  a baseline, which is another unsupervised learning-based  optimization method for problems over graphs.  We also compare our graph MaxCut results with the best known heuristic, BLS \cite{benlic2013breakout}. 

\subsection{Hypergraph MaxCut}
We used \sys{} to solve the hypergraph MaxCut problem. In  hypergraph MaxCut, the goal is to partition the nodes of a hypergraph into two groups such that the cut size is maximized, where the cut size is the number of hyperedges that have at least one pair of nodes belonging to  different groups. We compare the performance of \sys{} with simulated annealing and a direct gradient descent method (ADAM) as the baselines.  Also, note that PI-GNN \cite{schuetz2022combinatorial} can not solve problems over hypergraphs and thus, we do not compare with it. 
We perform our experiments on both synthetic and real hypergraphs.
The real hypergraphs are extracted from  the  American Physical Society (APS) \cite{APS_website}. See Supplementary Information Section V for more details on the real and synthetic hypergraphs used in our experiments. 
In Fig.\,3, we show the performance of \sys{} compared  to SA and ADAM. It can be seen that with almost the same performance, the run time of \sys{} is remarkably less than the other two methods. Notably, the run time of \sys{} grows linearly with the number of nodes, while the run time of SA grows quadratically.  

In Extended Data Figure\,1, we compare \sys{} with the bipartite GNN baseline. As depicted in the figure, for almost the same performance, \sys{} has a considerably reduced run time compared to the bipartite GNN method (up to 8 times improvement). The extended run time of the bipartite GNN can be attributed to the fact that the bipartite graph equivalent of a hypergraph with $N$ nodes and $K$ hyperedges would have $N+K$ nodes. Given that the  complexity of graph and hypergraph neural networks typically scale with the number of nodes, our experiments affirm that hypergraph neural networks are a more efficient alternative compared to bipartite graph neural networks.

\begin{figure}
    \centering
    \subfigure[Performance of \sys{} and SA (Synthetic Hypergraphs)]{\includegraphics[width=0.45\textwidth]{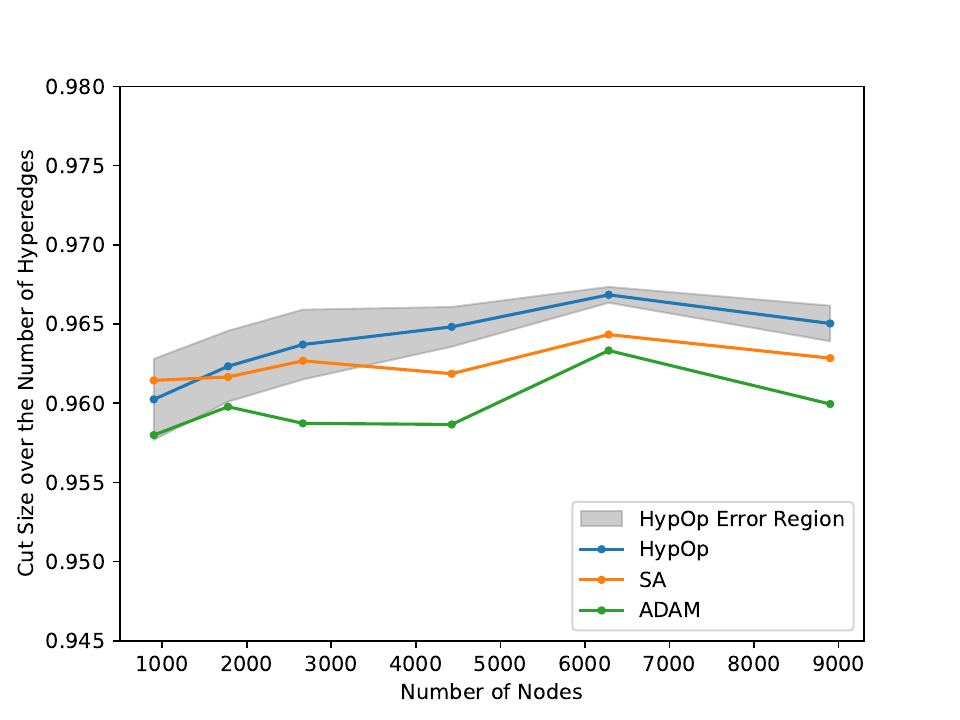}}
     \subfigure[Runtime of \sys{} and SA (Synthetic Hypergraphs)]{\includegraphics[width=0.45\textwidth]{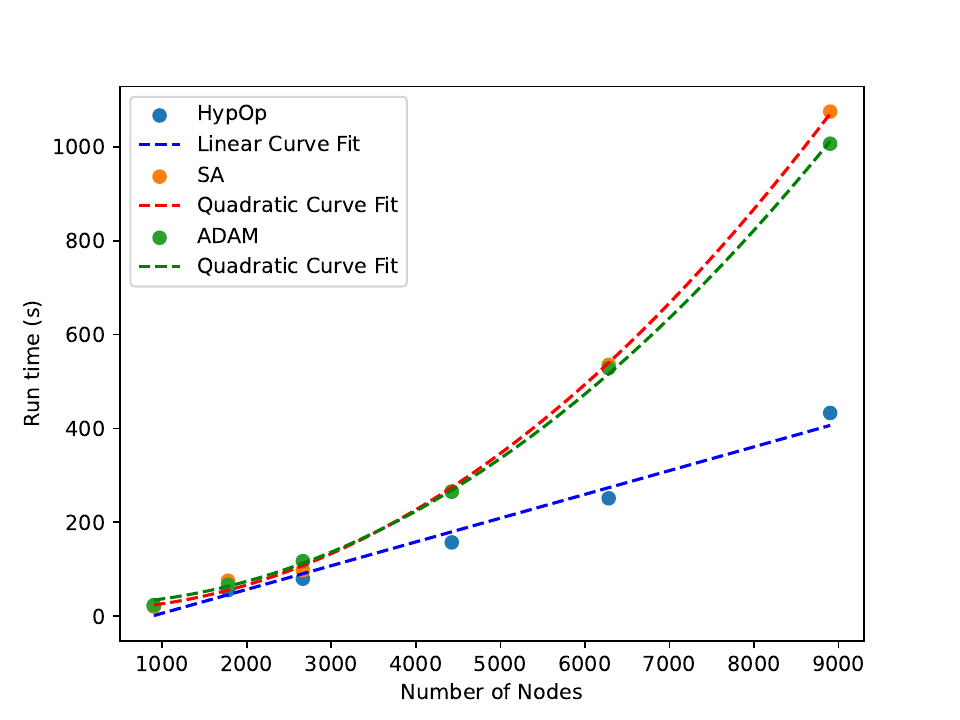}}\\
     \subfigure[Performance of \sys{} and SA (Real APS Hypergraphs)]{\includegraphics[width=0.45\textwidth]{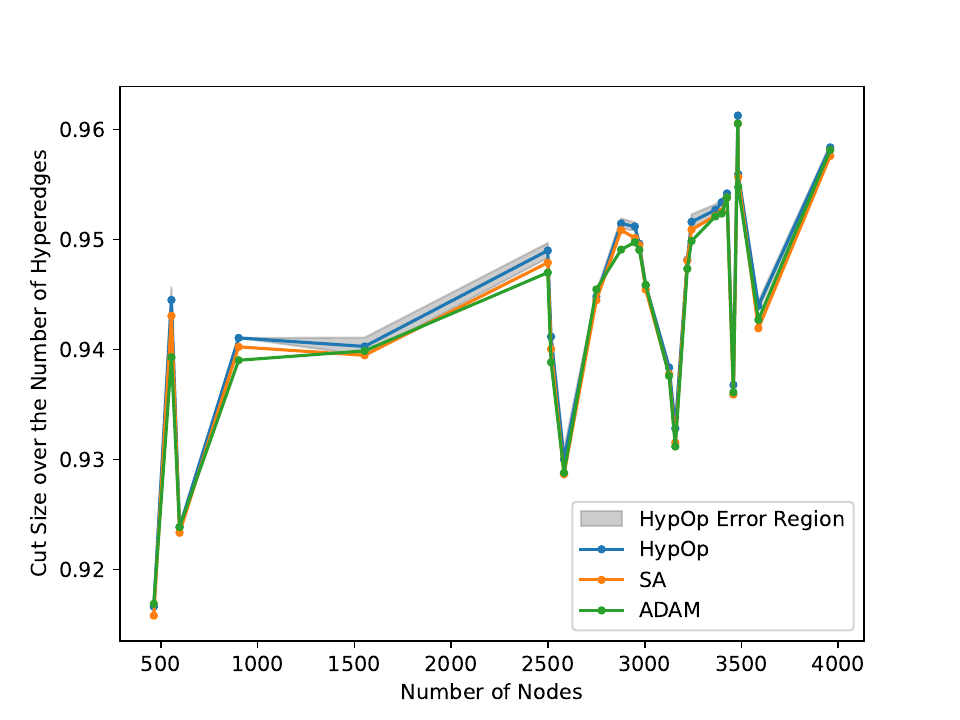}}
     \subfigure[Runtime of \sys{} and SA (Real APS Hypergraphs)]{\includegraphics[width=0.45\textwidth]{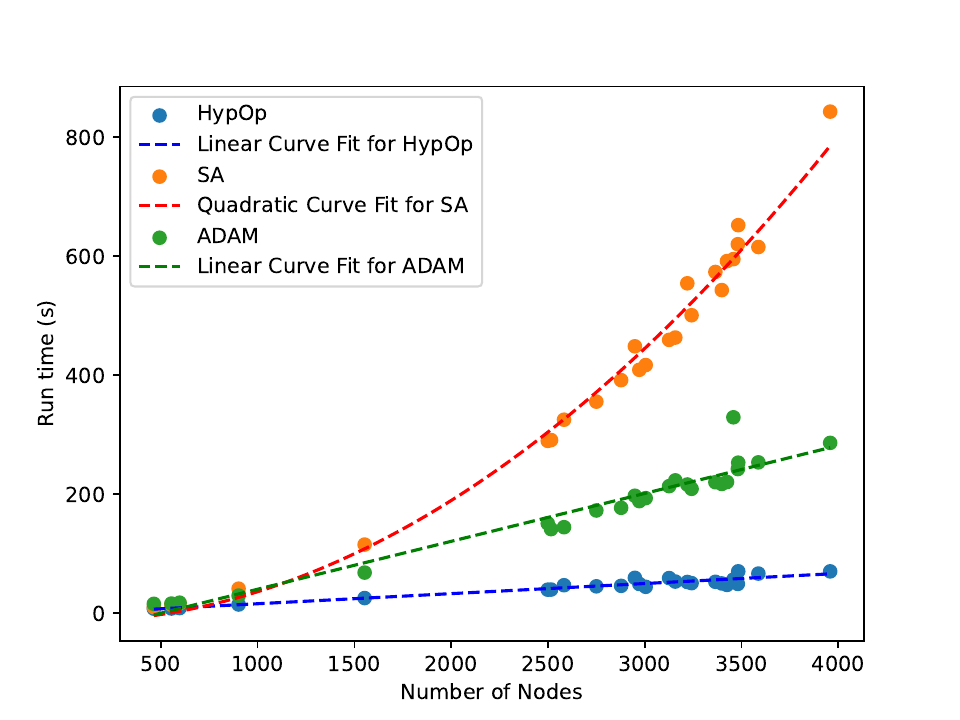}}
    \caption{Performance of \sys{} compared with SA and ADAM for hypergraph MaxCut problem on synthetic random hypergraphs and real APS hypergraphs. For almost the same performance, \sys{} has a remarkably less runtime compared to SA and ADAM, with SA's runtime growing quadratically fast, as opposed  to runtime of \sys{} that grows linearly. \sys{} performance is presented as the average of results from 10 sets of experiments and the error region  shows the standard deviation of the  results. The runtime standard deviations  were insignificant, with an average of 3.7s. }
    \label{fig:hyp_sa_syn}
\end{figure}

\subsection{Maximum Independent Set Problem}
In order to  compare \sys{} with the baseline work PI-GNN \cite{schuetz2022combinatorial} and  highlight the benefit of the fine-tuning step using SA in \sys{}, we solved  MIS problems over graphs using \sys{}.
In Fig.\,4, we show the performance and run time of PI-GNN and \sys{} over regular random graphs with $d=3$ and $d=5$. As can be seen in the figures, \sys{} performs  better than PI-GNN both in terms of optimality and run time. While generating better MIS sizes, the run time of \sys{} is up to five times better than PI-GNN. \sys{} also scales better than PI-GNN to larger graphs. Having different components in the optimization tool allows us to adjust the hyperparameters of them to better suit the specific problem at hand. Specially, we can have a more efficient algorithm by taking advantage of different types of solvers.
\begin{figure}
    \centering
\subfigure[MIS Ratio]{\label{Maxind_regular_d3-a}\includegraphics[width=0.47\textwidth]{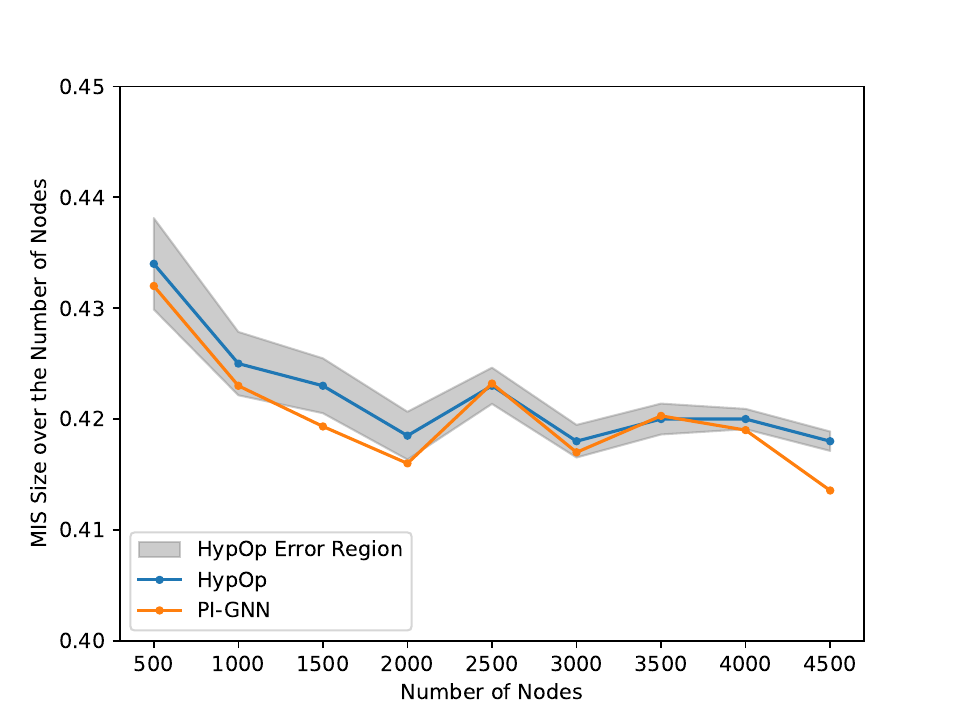}}
    \subfigure[Runtime]{\label{Maxind_regular_d3-b}\includegraphics[width=0.47\textwidth]{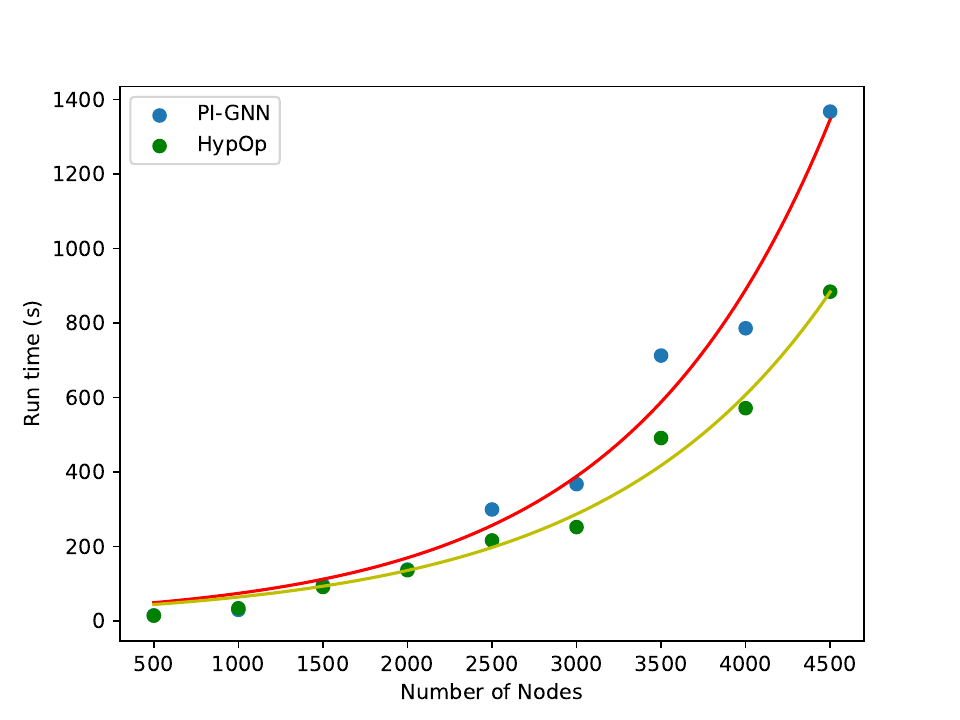}}\\
    \subfigure[MIS Ratio]{\label{Maxind_regular_d5-a}\includegraphics[width=0.47\textwidth]{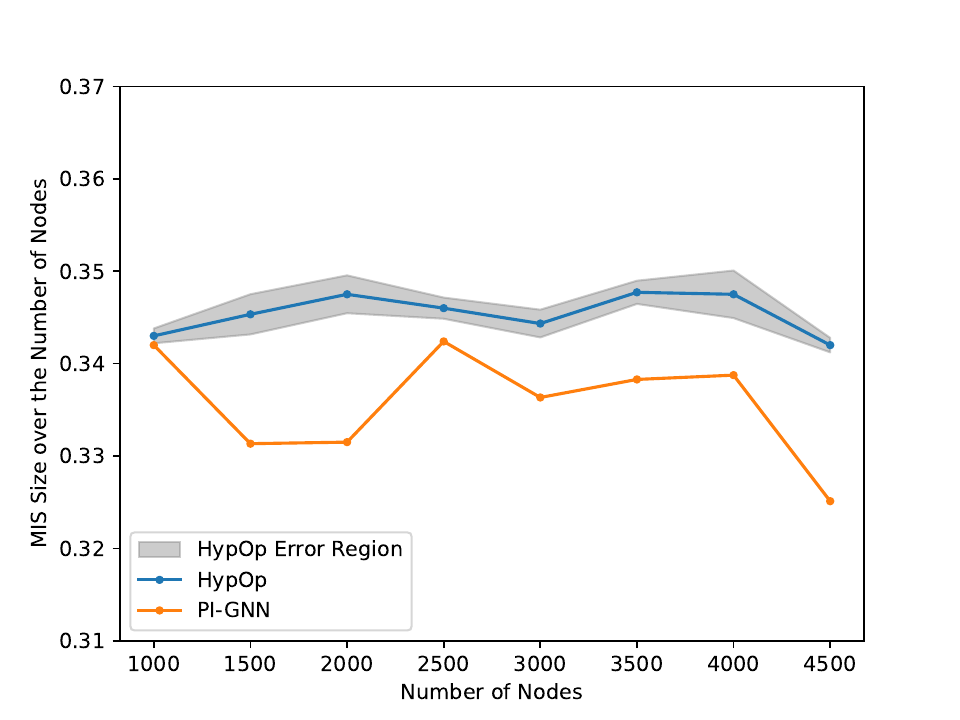}}
    \subfigure[Runtime]{\label{Maxind_regular_d3-b}\includegraphics[width=0.47\textwidth]{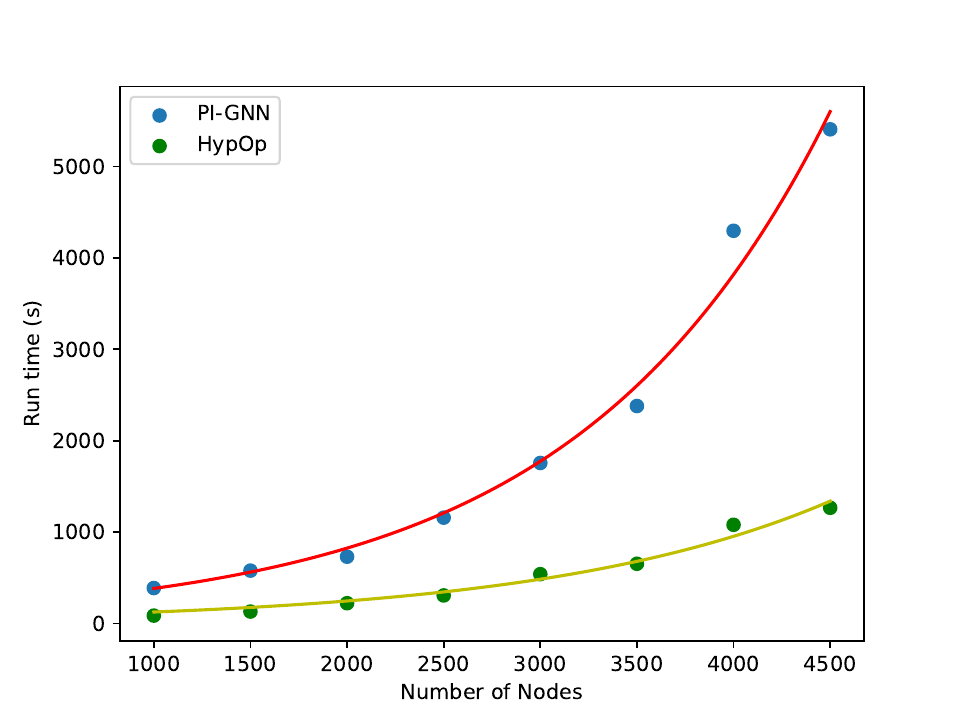}}
    \caption{Performance and runtime of \sys{} and PI-GNN for MIS problem on random regular graphs with $d=3$ and $d=5$. \sys{} achieves better results than PI-GNN in a  less amount of time.  \sys{} performance is presented as the average of results from 10 sets of experiments, with the error region  showing the standard deviation of the  results. The runtime standard deviations  were insignificant, with an average of 10s.  } 
    \label{fig:Maxind_regular_d3}
\end{figure}

\subsection{Parallel and Distributed Training}
In this section, we provide experimental results to  verify the effectiveness of the scalable multi-GPU HyperGNN training described in Algorithm\,\ref{alg:decentrailized_gpu} and Supplementary Information Algorithm\,2.  Note that all of the  prior  results were generated by single GPU training. For this section, we solved  MaxCut problem on one of the Stanford graphs, G22 \cite{gset}, in addition to  a larger graph OGB-Arxiv~\cite{hu2020open} with 169343 nodes and 1166243 edges. The results are summarized in Table~\ref{tab:multi-gpu-data}. For the smaller Stanford graph (G22), the single GPU training and Supplementary Information Algorithm\,2 show similar performance with improved run time using Algorithm\,2, while Algorithm\,\ref{alg:decentrailized_gpu} provides the outcome three times faster than the single GPU training. As the graph size becomes larger (OGB-Arxiv graph), single GPU training faces OOM (out of memory) error, whereas multi GPU training in \sys{} can solve this problem.  Algorithm\,\ref{alg:decentrailized_gpu} solves the problem almost two times faster than Algorithm\,2 with slight performance degradation ($\approx 0.1 \%$ degradation).

\begin{table}[]
\resizebox{\textwidth}{!}{%
\begin{tabular}{|c|c|c|cc|cc|cc|}
\hline
\multirow{2}{*}{Graph} & \multirow{2}{*}{Nodes} & \multirow{2}{*}{Edges} & \multicolumn{2}{c|}{Single GPU} & \multicolumn{2}{c|}{Parallel Multi-GPU} & \multicolumn{2}{c|}{Distributed   Multi-GPU} \\ \cline{4-9} 
 &  &  & \multicolumn{1}{c|}{Performance} & Time & \multicolumn{1}{c|}{Performance} & Time & \multicolumn{1}{c|}{Performance} & Time \\ \hline
Stanford G22 & 2,000 & 19,990 & \multicolumn{1}{c|}{13,185} & 15s & \multicolumn{1}{c|}{13,160} & 12s & \multicolumn{1}{c|}{13,128} & 5s \\ \hline
OGB-Arxiv & 169,343 & 1,166,243 & \multicolumn{1}{c|}{OOM} & - & \multicolumn{1}{c|}{855,505} & 750s & \multicolumn{1}{c|}{854,746} &  411s\\ \hline
\end{tabular}%
}
\caption{Multi-GPU Performance. Although parallel multi-GPU training can speed up the optimization and resolve the OOM error, the distributed training provides considerable speed up with its result being within 0.09\% of the  parallel Multi-GPU result. }
\label{tab:multi-gpu-data}
\end{table}

\section{Transfer Learning with \sys{} for Different Types of Problems}
\label{sec:transfer}
We explored the applicability of transfer learning within \sys{} to address various problems on a single graph. Our approach involves pre-training the HyperGNN on a specific problem type and then utilizing the pre-trained model for a different optimization problem on the same graph, focusing on optimizing the input embedding exclusively. In particular, we freeze the weight matrix $W$ and only optimize  the input embedding $\sigma$. The  idea is that the HyperGNN learns an effective function to capture graph features (modeled by $\mathcal{G}(\sigma,w)$) and therefore, it can potentially be used for other types of optimization problems on the same graph. The optimization on the embedding $\sigma$ yields customized solutions tailored to the specific optimization problem. To validate this, we trained \sys{} for the MaxCut problem and  transferred it to address the MIS problem. In Extended Data Figure\,2, we show the performance of \sys{} with and without transfer learning. The models are pre-trained for MaxCut problem and transferred to solve MIS problem on the same graph. As we can see in the figure, the transferred model can achieve comparable results in almost no time compared to the vanilla training (training from scratch). Similar results can be seen in Extended Data Figure\,3, where we test the transferability of \sys{}  from hypergraph MaxCut problem into hypergraph partitioning (MinCut) problem. The hypergraph partitioning or MinCut problem aims at cutting the hypergraph into two groups of nodes with almost the same size and the smallest possible cut size. Notably, transfer learning has even provided better results in some instances. 
These results further underscores the promising potential of utilizing unsupervised learning-based approaches to effectively and efficiently address complex optimization problems.

\section{Applications in Scientific Discovery}
\label{sec:science}
Combinatorial optimization is vital in diverse scientific and industrial applications, addressing challenges in resource allocation, experimental design, and decision-making processes. In this section, we show an example of how \sys{} can help with solving combinatorial optimization problems related to scientific discovery. We consider a hypothetical experimental design for national drug control, utilizing NDC dataset \cite{NDC, Benson-2018-simplicial} (see Supplementary Information Section VI for details). In this scenario, the goal is to select a set of substances for regulatory measures, dividing the dataset into two isolated groups of drugs. The objective is to determine the maximum number of substances requiring regulation to achieve this division. This problem can be formulated as a hypergraph MaxCut problem, which  we solved  using  \sys{} and SA. As shown in  Table \ref{tb:NDC}, \sys{} outperforms SA in terms of results and run time on the NDC hypergraph with 3767 nodes and 29810 hyperedges.

\begin{table}[]
    \centering
    \begin{tabular}{||c||c||c||}
    \hline & \sys{} & Simulated Annealing \\ \hline
      Cut Size  & 28849 $\pm$ 91 & 28472  $\pm$ 92\\ \hline
        Runtime (s) & 1679 &  3616\\ \hline
    \end{tabular}
    \caption{Performance of \sys{} on NDC dataset. Compared to SA, \sys{} can achieve better result in less amount of time.  }
    \label{tb:NDC}
\end{table}

\section{Discussion}
\label{sec:discussion}
The comparison among \sys{}, SA, ADAM, and PI-GNN in Section \ref{sec:exp} imparts several crucial insights. Firstly, relying solely on unsupervised learning-based methods for end-to-end optimization has demonstrated limitations, as observed in this work and \cite{angelini2022modern}. Combining such methods with generic optimization algorithms like SA substantially enhances performance, yielding acceptable solutions even when HyperGNNs or GNNs struggle to provide satisfactory results. Secondly, \sys{} and other unsupervised learning-based optimization methods can effectively serve as initializers for optimization variables in other solvers, such as SA. Our experiments demonstrate that if SA is not initialized with HyperGNN outputs, its run time performance notably suffers. Thirdly, our experiments confirm that HyperGNNs can function as effective learnable transformation functions, converting input embeddings into original optimization variables. This transformation acts as an acceleration for the solver, offering  more efficient optimization paths. Given the widespread use of gradient descent across various domains, such as in ML model training, utilizing HyperGNNs as transformation functions has the potential to enhance gradient descent-based optimization in numerous applications.


\section{Data Availability}
	In this paper, we  used publicly available datasets of APS \cite{APS_website}, NDC \cite{NDC}, Gset \cite{gset}, and SATLIB \cite{hoos2000satlib}, together with synthetic hypergraphs and graphs.  The procedure under which  the synthetic hypergraphs and graphs are generated is explained throughout the paper.  Some examples of the synthetic hypergraphs are provided with the code at \cite{Heydaribeni2024_ocean}. 

 \section{Code Availability}
 The code has been made publicly available at \cite{Heydaribeni2024_ocean}. We used Python 3.8 together with the following packages:  torch 2.1.1, tqdm 4.66.1, h5py 3.10.0, matplotlib 3.8.2, networkx 3.2.1, numpy 1.21.6, pandas 2.0.3, scipy 1.11.4, and sklearn 0.0. We used PyCharm 2023.1.2 and Visual Studio Code 1.83.1 softwares.

 \bibliographystyle{unsrt}

  \bibliography{Bibliography}

\begin{thebibliography}{10}

\bibitem{wang2023scientific}
Hanchen Wang, Tianfan Fu, Yuanqi Du, Wenhao Gao, Kexin Huang, Ziming Liu, Payal Chandak, Shengchao Liu, Peter Van~Katwyk, Andreea Deac, et~al.
\newblock Scientific discovery in the age of artificial intelligence.
\newblock {\em Nature}, 620(7972):47--60, 2023.

\bibitem{schuetz2022combinatorial}
Martin~JA Schuetz, J~Kyle Brubaker, and Helmut~G Katzgraber.
\newblock Combinatorial optimization with physics-inspired graph neural networks.
\newblock {\em Nature Machine Intelligence}, 4(4):367--377, 2022.

\bibitem{cappart2021combinatorial}
Quentin Cappart, Didier Ch{\'e}telat, Elias~B Khalil, Andrea Lodi, Christopher Morris, and Petar Veli{\v{c}}kovi{\'c}.
\newblock Combinatorial optimization and reasoning with graph neural networks.
\newblock {\em Journal of Machine Learning Research}, 24(130):1--61, 2023.

\bibitem{khalil2016learning}
Elias Khalil, Pierre Le~Bodic, Le~Song, George Nemhauser, and Bistra Dilkina.
\newblock Learning to branch in mixed integer programming.
\newblock In {\em Proceedings of the AAAI Conference on Artificial Intelligence}, volume~30, 2016.

\bibitem{bai2019simgnn}
Yunsheng Bai, Hao Ding, Song Bian, Ting Chen, Yizhou Sun, and Wei Wang.
\newblock Simgnn: A neural network approach to fast graph similarity computation.
\newblock In {\em Proceedings of the twelfth ACM international conference on web search and data mining}, pages 384--392, 2019.

\bibitem{gasse2019exact}
Maxime Gasse, Didier Ch{\'e}telat, Nicola Ferroni, Laurent Charlin, and Andrea Lodi.
\newblock Exact combinatorial optimization with graph convolutional neural networks.
\newblock {\em Advances in neural information processing systems}, 32:15580--15592, 2019.

\bibitem{nair2020solving}
Vinod Nair, Sergey Bartunov, Felix Gimeno, Ingrid Von~Glehn, Pawel Lichocki, Ivan Lobov, Brendan O'Donoghue, Nicolas Sonnerat, Christian Tjandraatmadja, Pengming Wang, et~al.
\newblock Solving mixed integer programs using neural networks.
\newblock {\em arXiv preprint arXiv:2012.13349}, 2020.

\bibitem{li2018combinatorial}
Zhuwen Li, Qifeng Chen, and Vladlen Koltun.
\newblock Combinatorial optimization with graph convolutional networks and guided tree search.
\newblock {\em Advances in neural information processing systems}, 31:537--546, 2018.

\bibitem{karalias2020erdos}
Nikolaos Karalias and Andreas Loukas.
\newblock Erdos goes neural: an unsupervised learning framework for combinatorial optimization on graphs.
\newblock {\em Advances in Neural Information Processing Systems}, 33:6659--6672, 2020.

\bibitem{toenshoff2021graph}
Jan Toenshoff, Martin Ritzert, Hinrikus Wolf, and Martin Grohe.
\newblock Graph neural networks for maximum constraint satisfaction.
\newblock {\em Frontiers in artificial intelligence}, 3:580607, 2021.

\bibitem{mirhoseini2021graph}
Azalia Mirhoseini, Anna Goldie, Mustafa Yazgan, Joe~Wenjie Jiang, Ebrahim Songhori, Shen Wang, Young-Joon Lee, Eric Johnson, Omkar Pathak, Azade Nazi, et~al.
\newblock A graph placement methodology for fast chip design.
\newblock {\em Nature}, 594(7862):207--212, 2021.

\bibitem{yolcu2019learning}
Emre Yolcu and Barnab{\'a}s P{\'o}czos.
\newblock Learning local search heuristics for boolean satisfiability.
\newblock {\em Advances in Neural Information Processing Systems}, 32:7992--8003, 2019.

\bibitem{ma2019combinatorial}
Qiang Ma, Suwen Ge, Danyang He, Darshan Thaker, and Iddo Drori.
\newblock Combinatorial optimization by graph pointer networks and hierarchical reinforcement learning.
\newblock {\em arXiv preprint arXiv:1911.04936}, 2019.

\bibitem{kool2018attention}
Wouter Kool, Herke van Hoof, and Max Welling.
\newblock Attention, learn to solve routing problems!
\newblock In {\em International Conference on Learning Representations}, 2018.

\bibitem{asghari2022transformation}
Mohammad Asghari, Amir~M Fathollahi-Fard, SMJ Mirzapour Al-E-Hashem, and Maxim~A Dulebenets.
\newblock Transformation and linearization techniques in optimization: A state-of-the-art survey.
\newblock {\em Mathematics}, 10(2):283, 2022.

\bibitem{feng2021hypergraph2}
Song Feng, Emily Heath, Brett Jefferson, Cliff Joslyn, Henry Kvinge, Hugh~D Mitchell, Brenda Praggastis, Amie~J Eisfeld, Amy~C Sims, Larissa~B Thackray, et~al.
\newblock Hypergraph models of biological networks to identify genes critical to pathogenic viral response.
\newblock {\em BMC bioinformatics}, 22(1):1--21, 2021.

\bibitem{murgas2022hypergraph}
Kevin~A Murgas, Emil Saucan, and Romeil Sandhu.
\newblock Hypergraph geometry reflects higher-order dynamics in protein interaction networks.
\newblock {\em Scientific Reports}, 12(1):20879, 2022.

\bibitem{zhu2018social}
Jianming Zhu, Junlei Zhu, Smita Ghosh, Weili Wu, and Jing Yuan.
\newblock Social influence maximization in hypergraph in social networks.
\newblock {\em IEEE Transactions on Network Science and Engineering}, 6(4):801--811, 2018.

\bibitem{xia2022novel}
Liqiao Xia, Pai Zheng, Xiao Huang, and Chao Liu.
\newblock A novel hypergraph convolution network-based approach for predicting the material removal rate in chemical mechanical planarization.
\newblock {\em Journal of Intelligent Manufacturing}, 33(8):2295--2306, 2022.

\bibitem{wen2012hyperspectral}
Yue Wen, Yue Gao, Shaohui Liu, Qimin Cheng, and Rongrong Ji.
\newblock Hyperspectral image classification with hypergraph modelling.
\newblock In {\em Proceedings of the 4th International Conference on Internet Multimedia Computing and Service}, pages 34--37, 2012.

\bibitem{feng2019hypergraph}
Yifan Feng, Haoxuan You, Zizhao Zhang, Rongrong Ji, and Yue Gao.
\newblock Hypergraph neural networks.
\newblock In {\em Proceedings of the AAAI conference on artificial intelligence}, volume~33, pages 3558--3565, 2019.

\bibitem{angelini2022modern}
Maria~Chiara Angelini and Federico Ricci-Tersenghi.
\newblock Modern graph neural networks do worse than classical greedy algorithms in solving combinatorial optimization problems like maximum independent set.
\newblock {\em Nature Machine Intelligence}, pages 1--3, 2022.

\bibitem{kirkpatrick1983optimization}
Scott Kirkpatrick, C~Daniel Gelatt~Jr, and Mario~P Vecchi.
\newblock Optimization by simulated annealing.
\newblock {\em science}, 220(4598):671--680, 1983.

\bibitem{kingma2014adam}
Diederik~P Kingma and Jimmy Ba.
\newblock Adam: A method for stochastic optimization.
\newblock {\em arXiv preprint arXiv:1412.6980}, 2014.

\bibitem{benlic2013breakout}
Una Benlic and Jin-Kao Hao.
\newblock Breakout local search for the max-cutproblem.
\newblock {\em Engineering Applications of Artificial Intelligence}, 26(3):1162--1173, 2013.

\bibitem{APS_website}
Aps dataset on {P}hysical {R}eview journals.
\newblock Published by the American Physical Society, https://journals.aps.org/datasets.

\bibitem{gset}
Y.~Ye.
\newblock The gset dataset.
\newblock https://web.stanford.edu/~yyye/ yyye/Gset/ (Stanford, 2003).

\bibitem{hu2020open}
Weihua Hu, Matthias Fey, Marinka Zitnik, Yuxiao Dong, Hongyu Ren, Bowen Liu, Michele Catasta, and Jure Leskovec.
\newblock Open graph benchmark: Datasets for machine learning on graphs.
\newblock {\em Advances in neural information processing systems}, 33:22118--22133, 2020.

\bibitem{NDC}
Ndc-substances dataset.
\newblock https://www.cs.cornell.edu/~arb/data/NDC-substances/.

\bibitem{Benson-2018-simplicial}
Austin~R. Benson, Rediet Abebe, Michael~T. Schaub, Ali Jadbabaie, and Jon Kleinberg.
\newblock Simplicial closure and higher-order link prediction.
\newblock {\em Proceedings of the National Academy of Sciences}, 2018.

\bibitem{hoos2000satlib}
Holger~H Hoos and Thomas St{\"u}tzle.
\newblock Satlib: An online resource for research on sat.
\newblock {\em Sat}, 2000:283--292, 2000.

\bibitem{Heydaribeni2024_ocean}
Nasimeh Heydaribeni, Xinrui Zhan, Ruisi Zhang, Tina Eliassi-Rad, and Farinaz Koushanfar.
\newblock (2024) {H}yp{O}p: Distributed constrained combinatorial optimization leveraging hypergraph neural networks [source code].
\newblock https://doi.org/10.24433/CO.4804643.v1.

\bibitem{bengio2021machine}
Yoshua Bengio, Andrea Lodi, and Antoine Prouvost.
\newblock Machine learning for combinatorial optimization: a methodological tour d’horizon.
\newblock {\em European Journal of Operational Research}, 290(2):405--421, 2021.

\bibitem{cummins2017end}
Chris Cummins, Pavlos Petoumenos, Zheng Wang, and Hugh Leather.
\newblock End-to-end deep learning of optimization heuristics.
\newblock In {\em 2017 26th International Conference on Parallel Architectures and Compilation Techniques (PACT)}, pages 219--232. IEEE, 2017.

\bibitem{liu2017deep}
Lu~Liu, Yu~Cheng, Lin Cai, Sheng Zhou, and Zhisheng Niu.
\newblock Deep learning based optimization in wireless network.
\newblock In {\em 2017 IEEE international conference on communications (ICC)}, pages 1--6. IEEE, 2017.

\bibitem{he2017deep}
Ying He, Zheng Zhang, F~Richard Yu, Nan Zhao, Hongxi Yin, Victor~CM Leung, and Yanhua Zhang.
\newblock Deep-reinforcement-learning-based optimization for cache-enabled opportunistic interference alignment wireless networks.
\newblock {\em IEEE Transactions on Vehicular Technology}, 66(11):10433--10445, 2017.

\bibitem{li2022overview}
Bingjie Li, Guohua Wu, Yongming He, Mingfeng Fan, and Witold Pedrycz.
\newblock An overview and experimental study of learning-based optimization algorithms for the vehicle routing problem.
\newblock {\em IEEE/CAA Journal of Automatica Sinica}, 9(7):1115--1138, 2022.

\bibitem{fawzi2022discovering}
Alhussein Fawzi, Matej Balog, Aja Huang, Thomas Hubert, Bernardino Romera-Paredes, Mohammadamin Barekatain, Alexander Novikov, Francisco~J R~Ruiz, Julian Schrittwieser, Grzegorz Swirszcz, et~al.
\newblock Discovering faster matrix multiplication algorithms with reinforcement learning.
\newblock {\em Nature}, 610(7930):47--53, 2022.

\bibitem{sun2022annealed}
Haoran Sun, Etash~Kumar Guha, and Hanjun Dai.
\newblock Annealed training for combinatorial optimization on graphs.
\newblock In {\em OPT 2022: Optimization for Machine Learning (NeurIPS 2022 Workshop)}, 2022.

\bibitem{lin2023comprehensive}
Haiyang Lin, Mingyu Yan, Xiaochun Ye, Dongrui Fan, Shirui Pan, Wenguang Chen, and Yuan Xie.
\newblock A comprehensive survey on distributed training of graph neural networks.
\newblock {\em Proceedings of the IEEE}, 2023.

\bibitem{chiang2019cluster}
Wei-Lin Chiang, Xuanqing Liu, Si~Si, Yang Li, Samy Bengio, and Cho-Jui Hsieh.
\newblock Cluster-gcn: An efficient algorithm for training deep and large graph convolutional networks.
\newblock In {\em Proceedings of the 25th ACM SIGKDD international conference on knowledge discovery \& data mining}, pages 257--266, 2019.

\bibitem{bai2021ripple}
Jiyang Bai, Yuxiang Ren, and Jiawei Zhang.
\newblock Ripple walk training: A subgraph-based training framework for large and deep graph neural network.
\newblock In {\em 2021 International Joint Conference on Neural Networks (IJCNN)}, pages 1--8. IEEE, 2021.

\bibitem{zeng2019graphsaint}
Hanqing Zeng, Hongkuan Zhou, Ajitesh Srivastava, Rajgopal Kannan, and Viktor Prasanna.
\newblock Graphsaint: Graph sampling based inductive learning method.
\newblock In {\em International Conference on Learning Representations}, 2019.

\bibitem{ramezani2021learn}
Morteza Ramezani and Weilin Cong.
\newblock Learn locally, correct globally: A distributed algorithm for training graph neural networks.
\newblock In {\em International Conference on Learning (ICLR)}, 2022.

\bibitem{javaheripi2020adans}
Mojan Javaheripi, Mohammad Samragh, Tara Javidi, and Farinaz Koushanfar.
\newblock Adans: Adaptive non-uniform sampling for automated design of compact dnns.
\newblock {\em IEEE Journal of Selected Topics in Signal Processing}, 14(4):750--764, 2020.

\bibitem{ding2020more}
Kaize Ding, Jianling Wang, Jundong Li, Dingcheng Li, and Huan Liu.
\newblock Be more with less: Hypergraph attention networks for inductive text classification.
\newblock In {\em Proceedings of the 2020 Conference on Empirical Methods in Natural Language Processing (EMNLP)}, pages 4927--4936, 2020.

\bibitem{larock2020understanding}
Timothy LaRock, Timothy Sakharov, Sahely Bhadra, and Tina Eliassi-Rad.
\newblock Understanding the limitations of network online learning.
\newblock {\em Applied Network Science}, 5:1--25, 2020.

\bibitem{rusch2023survey}
T~Konstantin Rusch, Michael~M Bronstein, and Siddhartha Mishra.
\newblock A survey on oversmoothing in graph neural networks.
\newblock {\em arXiv preprint arXiv:2303.10993}, 2023.

\bibitem{rong2019dropedge}
Yu~Rong, Wenbing Huang, Tingyang Xu, and Junzhou Huang.
\newblock Dropedge: Towards deep graph convolutional networks on node classification.
\newblock In {\em International Conference on Learning Representations}, 2019.

\end{thebibliography}

\newpage
\break 

\setcounter{figure}{0}
\renewcommand{\figurename}{Extended Data Figure}
\newcommand{\blue}[1]{\textcolor{blue}{ #1 }}

\begin{figure}
	\centering
	\subfigure[Performance]{\includegraphics[width=0.45\textwidth]{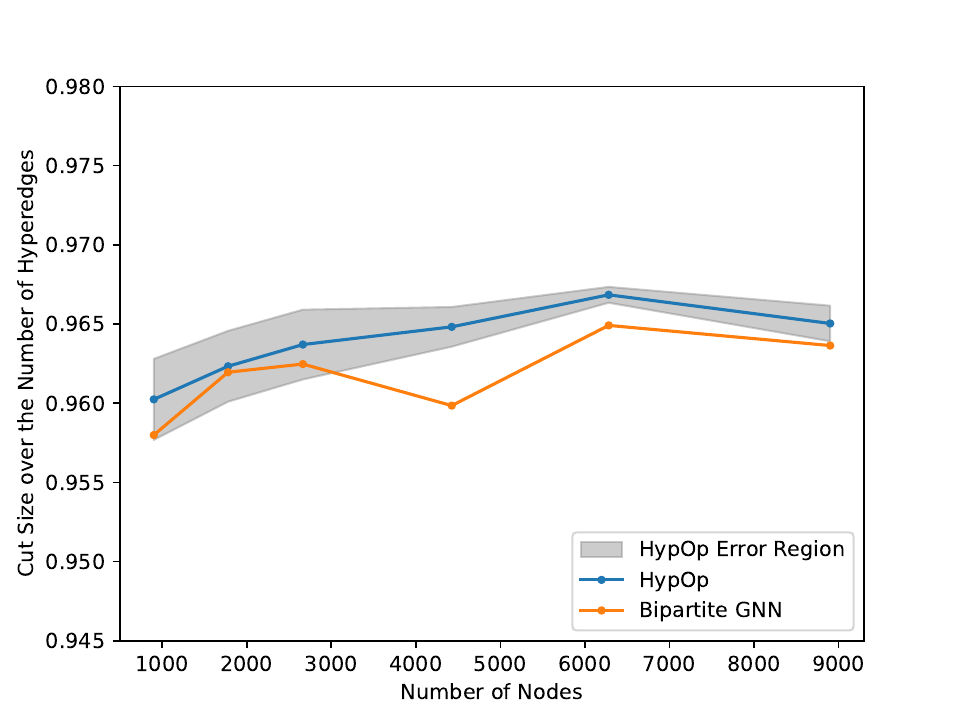}}
	\subfigure[Run time]{\includegraphics[width=0.45\textwidth]{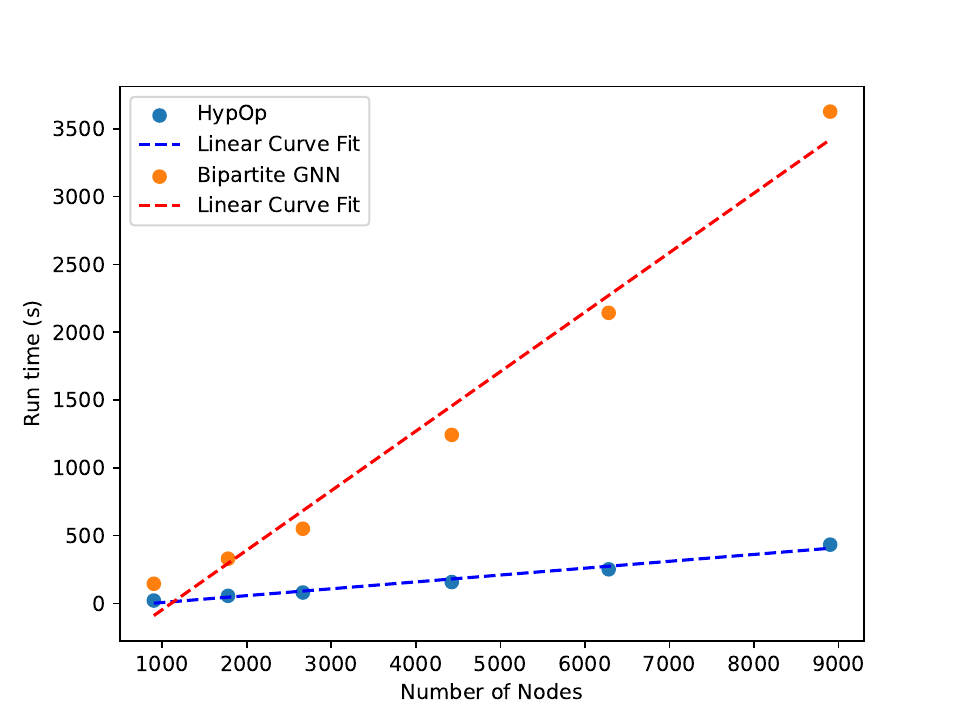}}
	\caption{\sys{} vs. Bipartite GNN. Comparison of \sys{} with the bipartite GNN baseline for hypergraph MaxCut problem on synthetic random  hypergraphs. For almost the same performance (a), \sys{} has a remarkably less run time compared to the bipartite GNN baseline (b). \sys{} performance is presented as the average of the results from 10 sets of experiments, with the error region  showing the standard deviation of the results.}
	\label{fig:hyp_bipartite}
\end{figure}

\begin{figure}
	\centering
	\subfigure[Performance]{\includegraphics[width=0.47\textwidth]{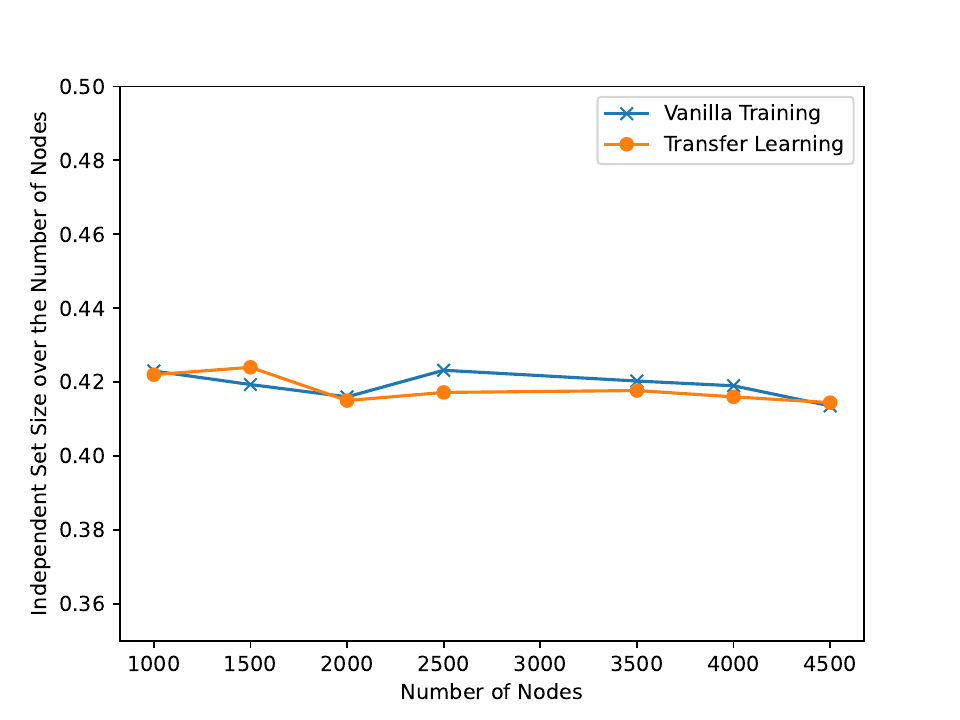}}
	\subfigure[Run time]{\includegraphics[width=0.47\textwidth]{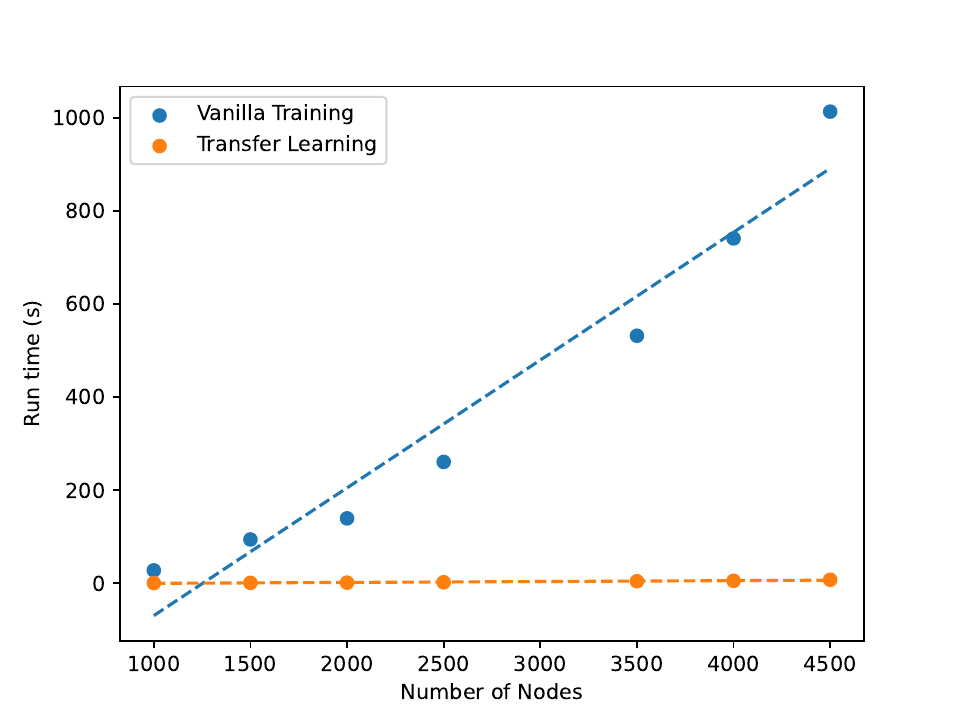}}
	\caption{Transfer Learning. Comparison of \sys{} with and without Transfer Learning  from MaxCut to MIS problem on random regular graphs with $d=3$.   For almost the same performance (a), transfer learning provides the results  in almost no amount of time compared to vanilla training (b). }
	\label{fig:transfer}
\end{figure}

\begin{figure}
	\centering
	\subfigure[Performance]{\includegraphics[width=0.47\textwidth]{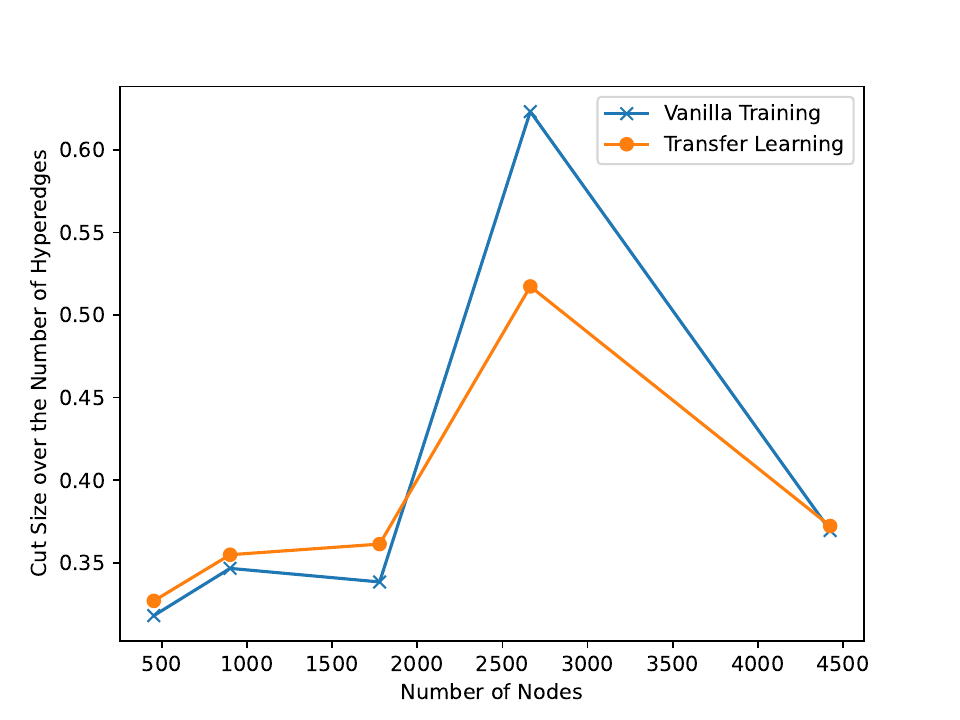}}
	\subfigure[Run time]{\includegraphics[width=0.47\textwidth]{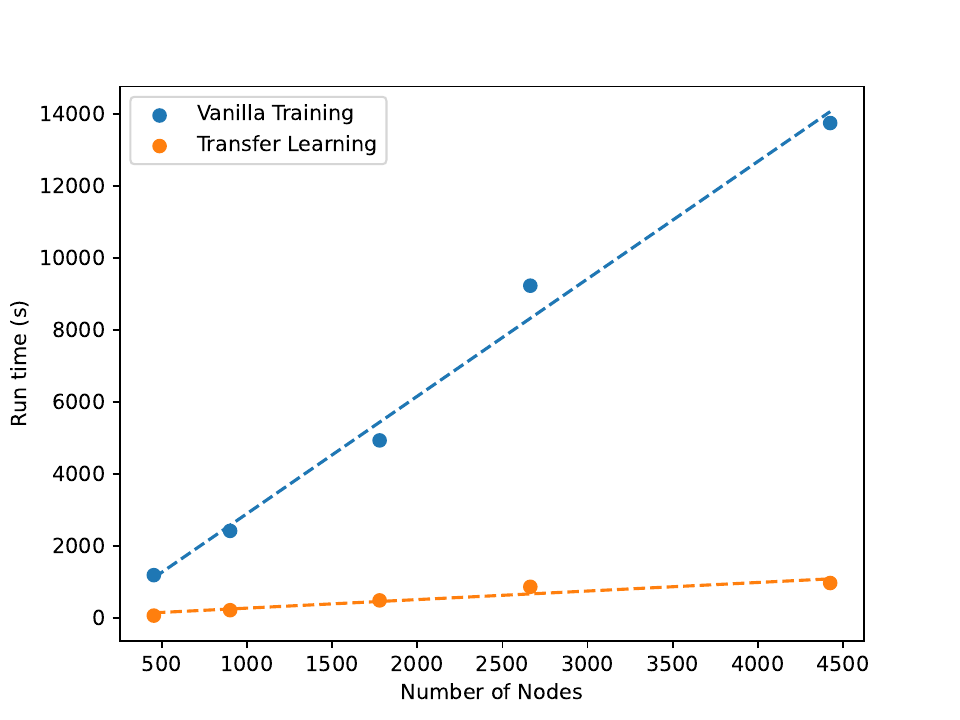}}
	\caption{Transfer Learning. Transfer Learning using \sys{} from Hypergraph MaxCut to Hypergraph MinCut on synthetic random hypergraphs. Compared to vanilla training, similar or better results are obtained using transfer learning (a) in a considerable less amount of time (b). Note that in the context of the Hypergraph MinCut problem, smaller cut sizes are favored. }
	\label{fig:transfer_hyper}
\end{figure}

    \newpage

 \addtocontents{toc}{\protect\fi}
   \clearpage

\title{Supplementary Information }


\begin{center}
    Supplementary Information
\end{center} 

\tableofcontents

\setcounter{figure}{0}
\renewcommand{\figurename}{Supplementary Figure}

\setcounter{table}{0}
\renewcommand{\tablename}{Supplementary Table}

\setcounter{page}{1}

\setcounter{section}{0}

\clearpage 
\section{Related Work}
The recent significant progress in AI and ML  has   prompted their applications in optimization \cite{ bengio2021machine, cummins2017end, liu2017deep, he2017deep}. Combinatorial optimization, given its intricate nature, stands to benefit greatly from learning-based approaches \cite{bengio2021machine, li2022overview, cappart2021combinatorial, mirhoseini2021graph, schuetz2022combinatorial}. Learning methods for combinatorial optimization problems can be categorized into supervised learning \cite{khalil2016learning, bai2019simgnn, gasse2019exact, nair2020solving, li2018combinatorial}, unsupervised learning \cite{schuetz2022combinatorial,karalias2020erdos, toenshoff2021graph}, and reinforcement learning \cite{mirhoseini2021graph, yolcu2019learning, ma2019combinatorial, kool2018attention}. Although supervised learning methods have shown promising results in various combinatorial optimization problems, they need a dataset of solved instances, posing challenges of computational intensity, expense, and limited generalizability to new or larger problems. 
Reinforcement learning is a powerful tool for optimization and has been proved to be effective for many problems. However, training RL models is usually  challenging and  time-consuming. Therefore, RL methods are best to be used for finding efficient algorithms to solve specific problems \cite{fawzi2022discovering}. 
Unsupervised learning methods enable fast and generalizable optimization of combinatorial optimization problems  without requiring a dataset of solved problem instances.  One issue with unsupervised learning methods is that they are based on gradient descent and therefore, they may converge to subpar local optima. 

To address the issues of the unsupervised learning methods getting stuck at local optima, \cite{sun2022annealed} proposed an annealed training framework for combinatorial optimization problems over graphs. A smoothing term is added to the loss function to help the training process escape the local optima. This smoothing strategy can be combined with  any unsupervised learning method for combinatorial optimization including the method described in this paper.    

While there is no existing research that models constrained combinatorial optimization problems as hypergraphs to capture higher-level dependencies between variables,  a study has approached maximum constraint satisfiability problems by modeling them as bipartite graphs and training a Graph Neural Network (GNN) to find solutions \cite{toenshoff2021graph}. However, it is limited to problems with constraints involving only two variables and is not suitable for problems with different types of constraints, such as 3-SAT.

\subsection{Distributed Training}
Numerous research works have addressed the acceleration and scalability of Graph Neural Network (GNN) training in distributed and parallel computing. Broadly categorized, these efforts focus on system-level optimizations, modifying graph data preprocessing and execution processes, or algorithm-level enhancements proposing efficient machine learning algorithms with minimal inter-server communication and improved model aggregation strategies \cite{lin2023comprehensive}. The system-level approaches are orthogonal to \sys{} and can be incorporated within \sys{} to further accelerate the training process. 
At the algorithm-level \cite{lin2023comprehensive}, works like \cite{chiang2019cluster, bai2021ripple, zeng2019graphsaint, ramezani2021learn} have focused on optimizing graph partitioning or model aggregation strategies for accelerated GNN training.
\sys{} doesn't enhance training by optimizing the graph partitioning scheme since it is designed for scenarios where the graph is inherently distributed across multiple servers and cases where the graph is partitioned beforehand. In this context, graph partitioning schemes are orthogonal to \sys{} and can be incorporated for improved scalability in the latter scenario. Furthermore, \sys{} exhibits two advantages over these approaches: (i) Unlike partition-based methods \cite{chiang2019cluster, bai2021ripple, zeng2019graphsaint}, \sys{} doesn't require subgraph sampling at every iteration, thereby reducing significant communication overhead; (ii) \cite{bai2021ripple} proposes a superior aggregation strategy to reduce communication overhead, however, it requires longer iterations for convergence. In contrast, \sys{} achieves comparable results within the same epochs as single GPU training, showcasing its efficiency and scalability for parallel/distributed training.
\begin{algorithm}[!ht]
\caption{Parallel Multi-GPU HyperGNN Training}
\label{alg:centralized_gpu}
\begin{algorithmic}[1]
\Require Hypergraph $G = (V, E)$, HyperGNN Model $M$, Number of GPUs $S$, Epoch Number $EP$
\Ensure Trained HyperGNN Model $M$

\State Initialize $M$ on each GPU
\For{$ep$ in $EP$}
    \State Shuffle $E$ to obtain $E_{\text{shuffled}}$
    \State Partition $E_{\text{shuffled}}$ into $S$ sublists: $E^1_{\text{shuffled}}, E^2_{\text{shuffled}}, \ldots, E^S_{\text{shuffled}}$
    
    \For{$s \leftarrow 1$ to $S$}
        \State Assign $E^s_{\text{shuffled}}$ to GPU$^s$
        \State Train $M^s$ on GPU$^s$ with edges $E^s_{\text{shuffled}}$
    \EndFor
    
    \State Aggregate gradients from all GPUs and update $M$
    \State Broadcast updated $M$ to all GPUs
\EndFor

\State \Return $M$

\end{algorithmic}
\end{algorithm}

\section{Preliminaries: Hypergraph  Neural Networks}

  \label{sec:prelim:hyp}
 A hypergraph is a generalized version of a graph in which, instead of having edges connecting two nodes, we have hyperedges connecting multiple nodes to each other,  allowing for more complex relationships than traditional graphs. Hypergraph neural network (HyperGNN) \cite{feng2019hypergraph}
 is an advanced variation of graph neural networks that is  designed to handle hypergraphs.  HyperGNNs extend the concept of GNNs to hypergraphs by introducing new aggregation mechanisms that consider higher-order connections (hyperedges). HyperGNNs can capture intricate dependencies among multiple nodes simultaneously, making them well-suited for tasks involving data with complex associations, such as social interactions, collaboration data, knowledge graphs, and recommendation systems. In hypergraph convolutional networks  introduced in \cite{feng2019hypergraph}, features of nodes are first aggregated in each hyperedge and are then passed on to new nodes for another aggregation. In particular, we have the following layer-wise operation:
 \begin{align}
    H^{(l+1)}=\sigma(D_v^{-\frac{1}{2}}{A}D_e^{-1}A^{T}D_v^{-\frac{1}{2}}H^{(l)}W^{(l)})
\end{align}
 where $H^{(l)}\in \mathbb{R}^{N\times D}$ is the matrix of node features at $l_{th}$ layer, $A$ is the hypergraph incidence matrix defined as $(A)_{ij}=1$ if node $i$ belongs to hyperedge $j$, and $(A)_{ij}=0$ otherwise. $D_v$ and $D_e$ are the diagonal degree matrices of nodes and hyperedges, respectively, and $W^{(l)}$ is the $l_{th}$ layer  trainable weight matrix. 
 
\section{HyperGCN Structure}
In  our HyperGCN model, we use two convolutional layers with a Relu activation function after the first layer. For the case of binary variables,  we  use Sigmoid activation function to generate a number between 0 and 1. The dimension of the convolutional layers are chosen as follows. The first layer has an $f$ dimensional input and outputs an $f/2$ dimensional vector. The second layer's input and output dimensions are $f/2$ and $1$, respectively. For most of out experiments, we set $f=\sqrt{N}$, where $N$ is the number of variables.
 
 \section{Black-Box Optimization of the Hyperparameters}
We have incorporated the optional hyperparameter optimization step in \sys{}, in which we conduct a sampling based black-box optimization w.r.t. the hyperparameters of the algorithm to improve the performance of \sys{}. We use AdaNS optimization tool  \cite{javaheripi2020adans} for this step. One of the hyperparameters that we optimize over is the initial value of the input embedding of the HyperGNN.  Since our method is based on gradient descent, the initialization point can potentially be important in avoiding bad local optima and therefore,  the output of the algorithm can be improved if we optimize w.r.t. the initial point. Other potential hyperparameters that we can optimize with this tool  are the learning rate and the SA hyperparameters. 

\section{Real and Synthetic Hypergraphs}
The real hypergraphs  used in our experiments are extracted from the the  American Physical Society (APS) \cite{APS_website}. The  APS dataset  contains the basic metadata of all APS journal articles from year 1993 to 2021. Focusing on a specific journal within this dataset, Physical Review E (PRE), we partitioned the data into 3-year segments, and extracted the  giant connected component of the authorship hypergraphs. The synthetic hypergraphs are generated randomly for  specific values of number of nodes, edges, and upper bounds on the node degrees.

\section{National Drug Code (NDC) Dataset}
NDC dataset consists of the substances that make up drugs in the National Drug
Code Directory maintained by the Food and Drug Administration \cite{NDC, Benson-2018-simplicial}. This dataset consists of a list of drugs and the substances that make them. The drugs are considered as nodes and  substances are presented by hyperedges. Drugs (nodes) within each hyperedge contain that substance. We considered the experiment of choosing a set of substances such that if one puts some regulatory measures on  them, the rest of the dataset is divided into two separate (smaller) isolated group of drugs that can be addressed separately.  We are interested to know the maximum number of substances that we need to regulate to create two separate datasets. This experiment can be represented  by a hypergraph MaxCut problem. Before applying \sys{} on NDC dataset,  We first removed the hyperedges containing one node and then removed the isolated nodes. As a result, we constructed a hypergraph with 3767 nodes and 29810 hyperedges. 

\section{Experiment Results}
\subsection{Graph MaxCut}
We have  considered graph MaxCut problem to compare \sys{} with the baseline solver \cite{schuetz2022combinatorial} (PI-GNN). Similar to PI-GNN, we run \sys{} on Gset dataset \cite{gset}.  Table \ref{tb:maxcut} summarizes our experiments on this dataset.
As can be seen in the table, \sys{} outperforms PI-GNN in all of the tested graphs.

We need to note here that Graph MaxCut is a very well studied problem and efficient heuristic algorithms have been developed to solve this problem. Although one can use algorithms such as \sys{} and PI-GNN to solve MaxCut problems, these algorithms are not expected to outperform the efficient heuristics that are customized to the  MaxCut problem specifics. However, the advantage of the learning-based algorithms such  as \sys{} and PI-GNN is that they are general solvers and one can use the same tool to solve thousands of optimization problems without the need to learn to run specific heuristics for each problem. So the main point of benchmarking \sys{} on MaxCut problem is to compare its performance with the baseline algorithm PI-GNN and to further motivate the additional components that we have added to the baseline method.

\begin{table}[h!]
\centering
\begin{normalsize}
\begin{tabular}{ ||c|| c|| c|| c|| c|| c || c|| c|| c||}
 \hline Graph & G14 & G15 & G22 & G49 & G50 & G55 \\ 
\hline Nodes  & 800  & 800  & 2000  & 3000  & 3000   & 5000    \\ 
\hline Edges & 4694 & 4661 & 19990 &  6000 & 6000  & 12498 \\
\hline Run-CSP & 2943 & 2928 & 13028 & 6000 & 5880 & 10116 \\
\hline PI-GNN &  3026  & 2990  & 13181 & 5918 & 5820  &10,138\\ 
\hline \sys{} &  3042  & 3030  & 13185 & 5999 & 5879 & 10150 \\
\hline Best Known Result & 3064 & 3050 & 13359 & 6000 & 5880 & 10,294  \\
\hline
\end{tabular}
\end{normalsize}  
\caption{Graph MaxCut with \sys{}. Performance of \sys{} for graph MaxCut problem on the Gset dataset. \sys{} outperforms other unsupervised learning based methods, while achieving results that are close to the best known heuristics. }
\label{tb:maxcut}
\end{table}

\subsection{SAT Problems}
We have tested our algorithm on random 3-SAT benchmarks from SATLIB dataset \cite{hoos2000satlib}. Table \ref{tb:sat} in the supplementary  section shows the performance of \sys{} on some of the satisfiable instances. We see that except for one with less than $2\%$ violation of constraints,  all others were successfully solved. We provide our experiments on SAT problems to showcase the success of our algorithm on various  problem instances, including the  challenging and widely studied SAT problems. We need to note that we do not expect \sys{} to compete with the powerful heuristics that are designed specifically for SAT problems and therefore, we do not provide any baselines in this section.
\begin{table}[h!]
\centering
\begin{normalsize}
\begin{tabular}{ ||c|| c|| c|| c|| c|| }
 \hline  SAT Dataset & uf 20-91 & uf 100-430 & uf 200-860 & uf 250-1065 \\ \hline 
Number of problems & 100 & 100 & 100 & 100 \\ \hline  
Percentage Solved  & $100\%$  & $100\%$ & $100\%$ & $98\%$ \\
\hline  Average number of unsatisfied clauses & 0 & 0 & 0 & 1 \\
\hline  Average time & 4s & 15s & 37s & 60s \\ 
\hline  Median time & 4s & 14s & 21s & 29s \\ \hline
\end{tabular}
\end{normalsize}
\caption{SAT with \sys{}. Performance of \sys{} for random 3-SAT problems in the SATLIB dataset.}
\label{tb:sat}
\end{table}

 \subsection{Resource Allocation Problem}
We use \sys{} to solve a resource allocation problem where there are a number of tasks to be assigned to a number of agents with given energy and budget constraints. In particular, each agent has an energy budget to spend on the tasks and each task requires a minimum energy to be completed. We can also consider a global arbitrary objective function to minimize while satisfying the constraints. One can model this problem as a hypergraph discovery problem where the hypergraph nodes correspond to the agents and the hyperedges correspond to the tasks. The nodes (agents) within each hyperedge (task) are assigned to that task. One needs to  discover different variations of this hypergraph to find feasible solutions and optimize the objective function. For simplicity, we do not consider a global objective function in our experiments and only try to find feasible solutions. 
In order to get realistic energy and budget constraints, we consider the APS co-authorship hypergraphs \cite{APS_website}. The papers correspond to tasks and the authors correspond to agents. The energy requirement of each task (paper) is the number of authors in it. Each agent has a budget that is a fixed number more than the papers that the agent appears in. We refer to this number as the budget surplus. The  feasible set would be larger if we increase the budget surplus. We use the same budget surplus for all of the hypergraphs that we have considered so that  as the number of agents in the hypergraph increases, the size of the feasible set shrinks.  We use the initial co-authorship hypergraphs to build the HyperGNN in \sys{} and try to find feasible solutions by assigning agents to tasks. We compare the result of \sys{} with SA in Supplementary Figure \ref{fig:task}. As can be seen in the figure, SA fails to find feasible solutions and its performance gets worse as the size of the feasible set shrinks (the number of agents increases). However, \sys{} can find feasible solutions for all number of agents. Furthermore, in terms of run time, \sys{} can solve these resource allocation problems in a significantly less time than SA (See Table \ref{tab:task}). 
\begin{table}[h]
\begin{minipage}[c]{0.53\textwidth}%
    \includegraphics[width=1\textwidth]{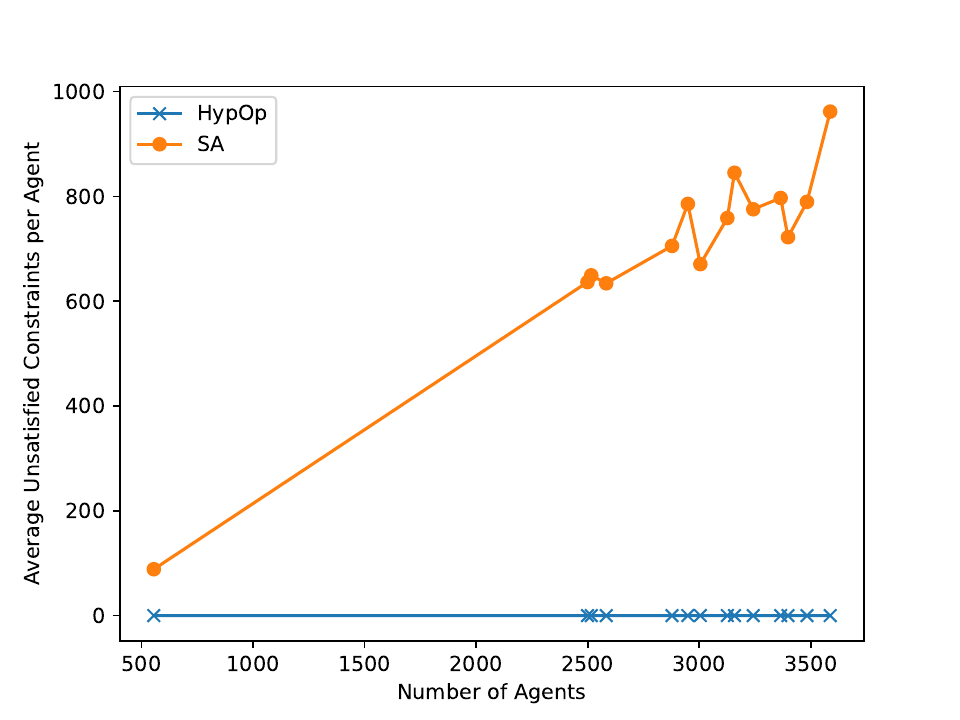} 
\captionof{figure}{Resource Allocation with \sys{}. Unsatisfied constraints in the Resource Allocation problem with different number of agents. \sys{} can find feasible solutions (zero constraint violation), while SA fails to do so.}
\label{fig:task}
\end{minipage}
 \begin{minipage}[c]{0.45\textwidth} 
  \caption{Average run time of \sys{} for resource allocation problems.}
 \label{tab:task}
\centering
\begin{normalsize}
    \begin{tabular}{|c|c|}
           \hline &  Run Time (s)   \\ \hline 
           \sys{} & 4668 $\pm$ 2017 \\ \hline 
           SA & 35472 $\pm$ 61448 \\  \hline
        \end{tabular}
\end{normalsize}
\end{minipage}
\end{table}


\section{HyperGNN Architecture Choice}
Two prominent HyperGNN architectures are Hypergraph Convolutional Networks (HyperGCNs) and Hypergraph Attention Networks (HyperGATs). HyperGCNs extend traditional graph convolutional networks to hypergraphs, enabling the incorporation of higher-order relationships.  On the other hand, HyperGATs leverage attention mechanisms to weigh the importance of different hyperedges and nodes, allowing the network to focus on relevant information during message passing.
In \sys{}, we opted for  HyperGCN architecture due to  the unsupervised nature of our approach and the fact that they   facilitate easier training for each problem instance.  To support our choice, we also implemented and tested a HyperGAT architecture within \sys{}. The HyperGAT model is based on the implementation in \cite{ding2020more}. In Supplementary Figure \ref{fig:attention}, we show the performance and runtime of \sys{} with HyperGCN against the same method with a HyperGAT architecture. As depicted in the figure, HyperGAT can not achieve the same performance compared to HyperGCN, while requiring  significantly more time to train.
\begin{figure}[t]
    \centering
    \subfigure[Performance]{\includegraphics[width=0.45\textwidth]{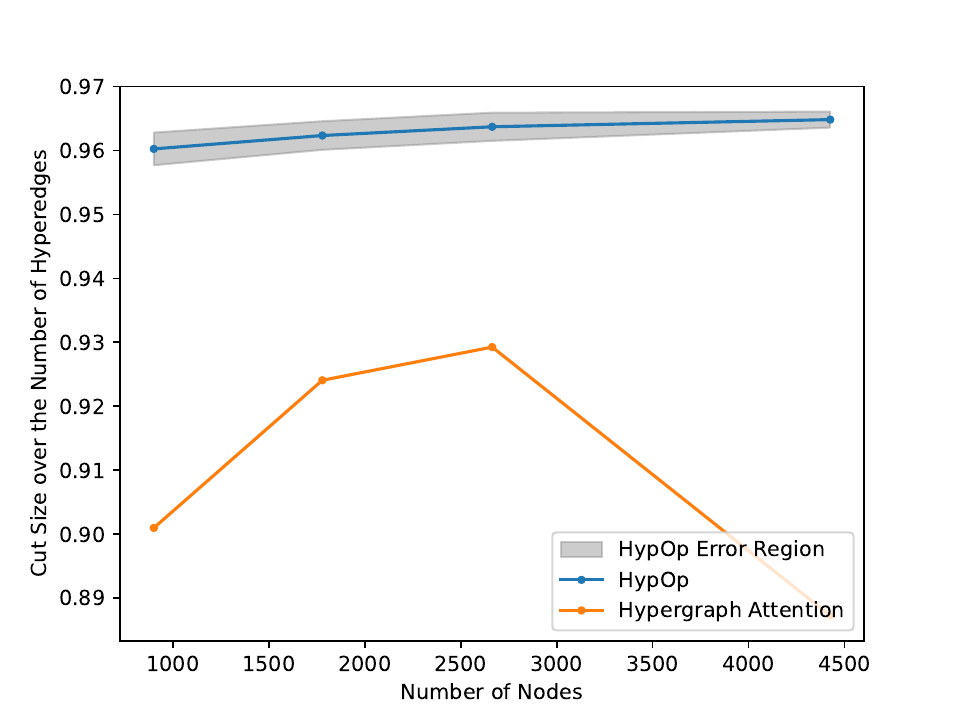}}
     \subfigure[Run time]{\includegraphics[width=0.45\textwidth]{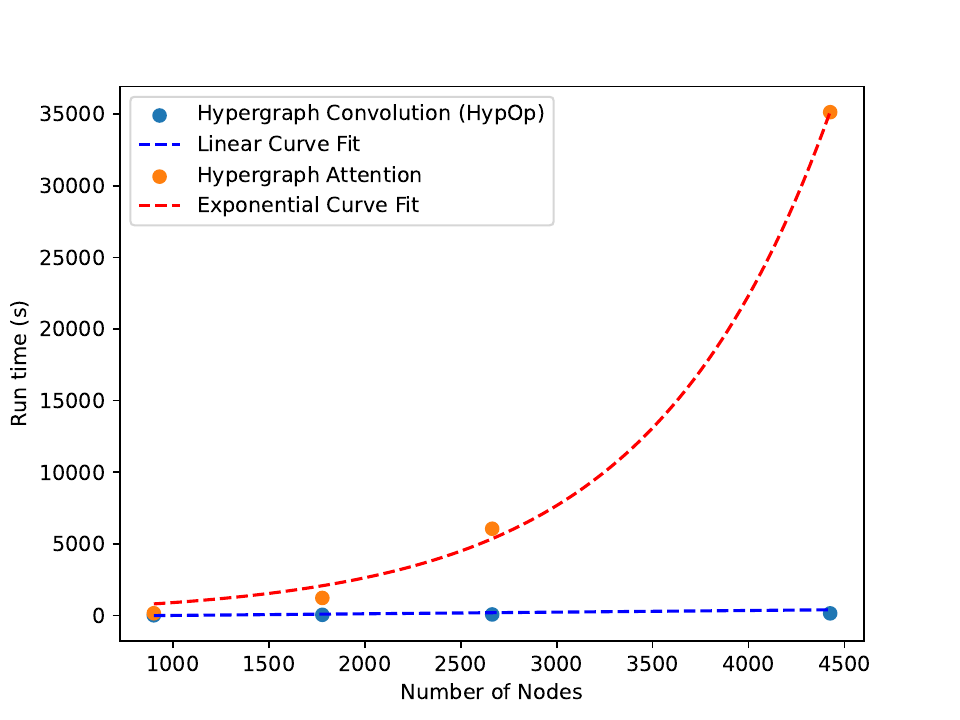}}
     \caption{HyperGCN vs. HyperGAT. Performance of \sys{} with HyperGCN architecture, compared with HyperGAT  for hypergraph MaxCut problem on synthetic random hypergraphs. HyperGAT can not achieve the  same performance compared to  HyperGCN, while requiring significantly more time to train.  \sys{} performance is presented as the average of the results from 10 sets of experiments and   the error region  shows the standard deviation of the results. }
    \label{fig:attention}
     \end{figure}

\section{Limitations: Phase Transition in GNN Performance}
\label{sec:limitations}
Although learning-based approaches for optimization have shown tremendous potential, one  has to be mindful of their limitations. 
In this section, we provide our experimental  results on the performance of GNNs for MIS problems on various graphs. Note that for the analysis in this section, we  consider the GNN outcome in \sys{} with a simple threshold mapping as used in PI-GNN. For problems on graphs (including MIS), this will be the same as the PI-GNN result. We compare this result to the fine-tuned outcome generated by \sys{} as an end-to-end solver. Our goal is to show the limitations of PI-GNN on some graphs. This analysis is inspired by the results of \cite{larock2020understanding}, where the authors study the limitations of the network online learning algorithms by considering different types of graphs.  


 We investigated different structures of graphs, namely, Erdos-Renyi Random graphs, Regular graphs, and Powerlaw graphs (graphs with powerlaw degree distribution), with various options for their number of nodes and edges. We observed that for MIS problem, PI-GNN can not learn anything on the graphs that are denser than a given threshold. In Supplementary Figure \ref{fig:phase:r:a} and \ref{fig:phase:r:c}, we show the performance of PI-GNN on Erdos-Renyi random graphs for different parameters, $p$, and different number of nodes, $N$, for two different learning rates. As is evident in the figure, the performance drops to 0 at a given threshold. The threshold, however, varies with the number of nodes. In Supplementary Figure \ref{fig:phase:r:b} and \ref{fig:phase:r:d}, we show the plot of the phase transition threshold w.r.t. the number of nodes. The Erdos-Renyi graphs are  well known for their phase transition in their connectivity. The threshold for this phase transition is $p^*=\frac{\text{ln}(N)}{N}$. As we have shown in Supplementary Figure \ref{fig:phase:r:b} and  \ref{fig:phase:r:d}, the phase transition threshold in the performance of PI-GNN is almost the same as the threshold for connectivity phase transition. This means that PI-GNN does not learn a solution over graphs that are dense enough to be connected. We provide our experimental results for two different learning rates to show that the performance drop is not caused by a high learning rate and is instead caused by the potential inability of GNNs to learn on dense graphs. 
\begin{figure}[h!]
    \centering
    \subfigure[Phase transition in the performance of PI-GNN over random graphs.]{ \label{fig:phase:r:a}  \includegraphics[width=0.45\textwidth]{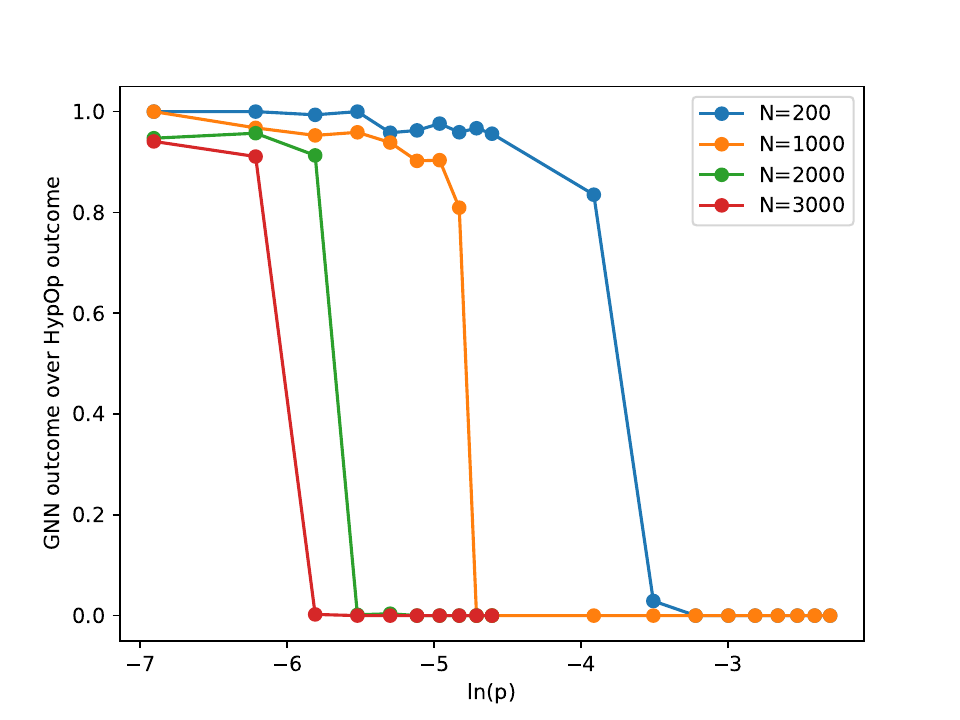}}
\subfigure[Phase transition threshold w.r.t. the graph size, and comparison with the Erdos-Renyi graph connectivity phase transition, $p^*=\frac{ln(N)}{N}$. ]{\label{fig:phase:r:b} \includegraphics[width=0.45\textwidth]{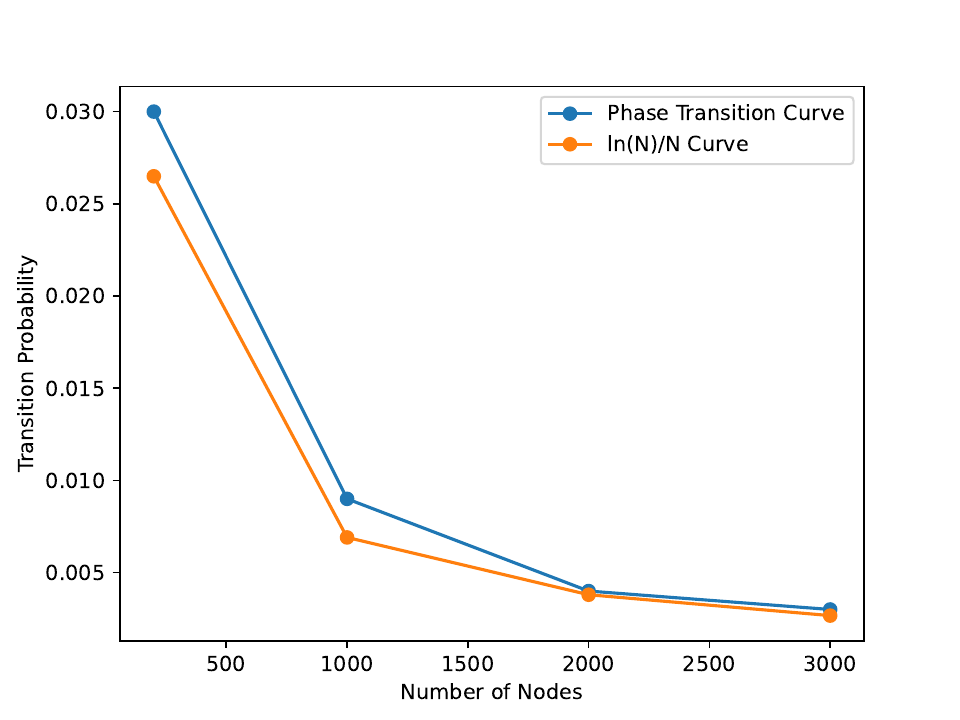}}
    \subfigure[Phase transition in the performance of PI-GNN over random graphs.]{ \label{fig:phase:r:c}  \includegraphics[width=0.45\textwidth]{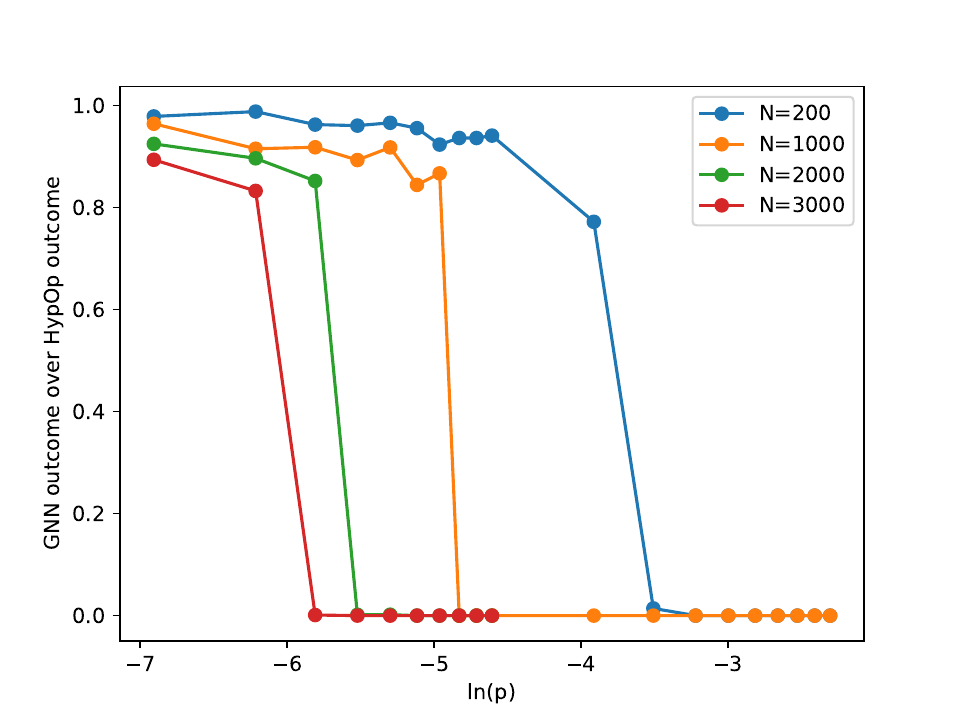}}
\subfigure[Phase transition threshold w.r.t. the graph size, and comparison with the Erdos-Renyi graph connectivity phase transition, $p^*=\frac{ln(N)}{N}$. ]{\label{fig:phase:r:d} \includegraphics[width=0.45\textwidth]{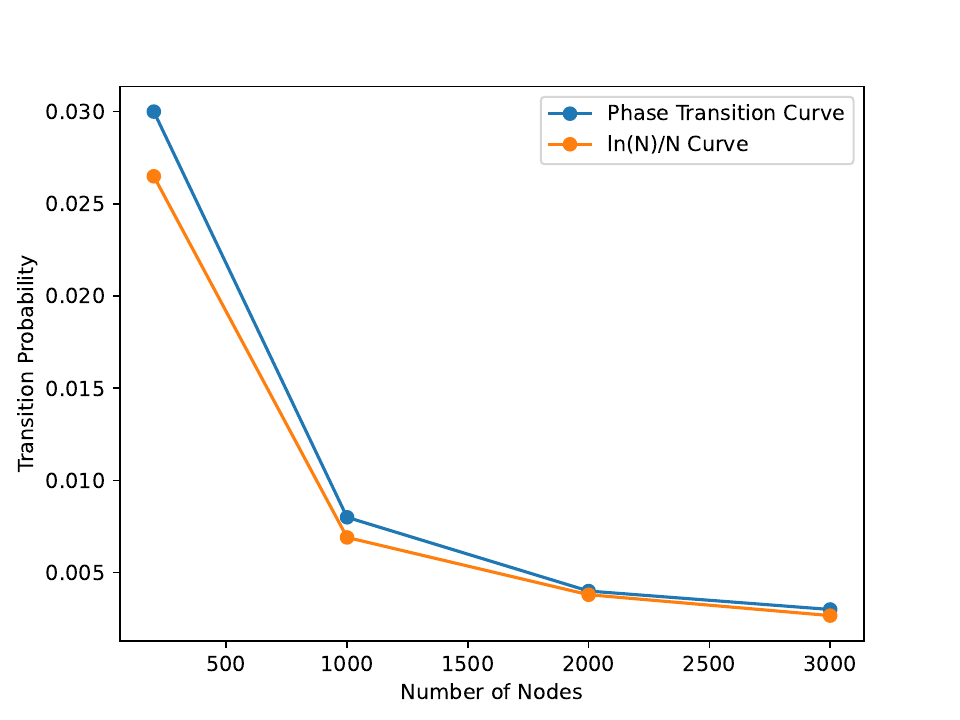}}
    \caption{Phase Transition in GNN Learning. Phase transition in the GNN-based optimization performance for MIS problem with learning rate of $1e-4$ in (a) and (b) and learning rate of $1e-5$ in (c) and (d).}
    \label{fig:phase:r}
\end{figure}

In order to investigate whether or not the specific structure of the graphs plays a role in the observed phase transition, we have compared the performance of PI-GNN on  Powerlaw graphs, Regular graphs, and  Random (Erdos-Renyi) graphs. In order to do the comparison, for each parameter of the powerlaw graphs, we  construct random and regular  graphs with almost the same number of edge density and compare the performance of PI-GNN over the three generated graphs. In Supplementary Figure \ref{fig:phase:st}, we show that the phase transition is the same for all of the three types of graphs when the performance is plotted w.r.t. the density ($p$) of the Random graph. Therefore, based on the existing evidence, we conclude that the performance is affected mainly by the density of the graph as opposed to its structure. We need to note here that most of the real graphs are powerlaw graphs with components between 2 to 4. Such graphs usually have lower densities than the observed critical density threshold for GNN learning. 
\begin{figure}[h!]
    \centering
\includegraphics[width=0.47\textwidth]{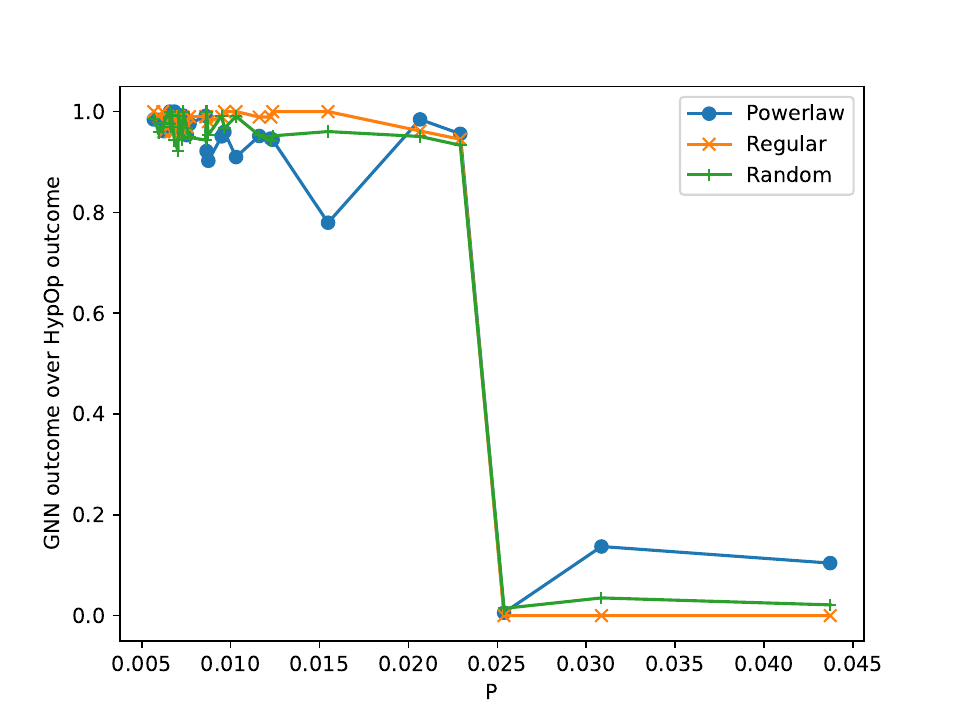}
    \caption{Phase Transition in GNN Learning. Phase transition in the performance of PI-GNN for Powerlaw, Regular, and Erdos-Renyi graphs.}
    \label{fig:phase:st}
\end{figure}

The phase transition in the GNN learning can be attributed to the oversmoothing issue that is widely observed in the GNNs throughout the literature \cite{rusch2023survey}. Roughly speaking, oversmoothing occurs due to too many  feature averaging, leading all hidden features to converge to the same point. In Supplementary Figure \ref{fig:smooth}, we show the maximum distance between the output embeddings of the different nodes after the main operations of the HyperGNN (two convolutional and two aggregation operators). We provide two saparate figures, one for the graphs for which the PI-GNN is successful in producing good outcomes (graphs that have density less than the threshold found in Supplementary Figure \ref{fig:phase:r}), and one for the graphs for which the PI-GNN is not successful in producing good outcomes (graphs that are denser than the threshold found in Supplementary Figure \ref{fig:phase:r}). We see that the embedding distance is decreasing and converges to a small number for the dense graphs. This observation confirms the occurrence of oversmoothing in the HyperGNN training for dense graphs. 

\begin{figure}[h!]
    \centering
    \subfigure[Successful GNN training: graphs that are sparser than the phase transition threshold.]{ \label{fig:smooth:a}  \includegraphics[width=0.45\textwidth]{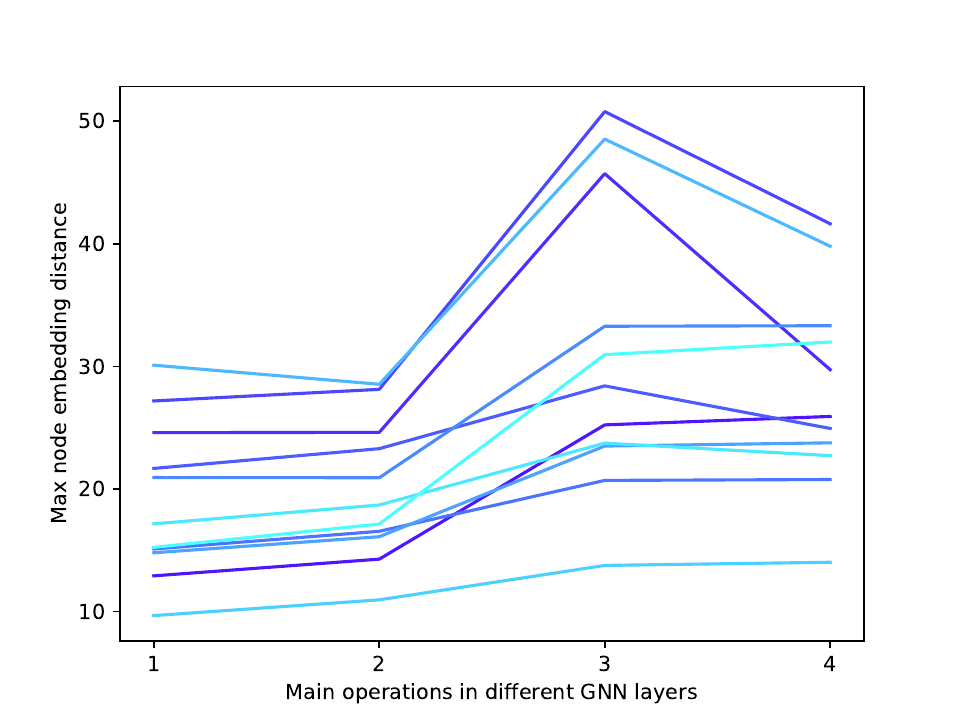}}
\subfigure[Unsuccessful GNN training: graphs that are denser than the phase transition threshold.]{\label{fig:smooth:b} \includegraphics[width=0.45\textwidth]{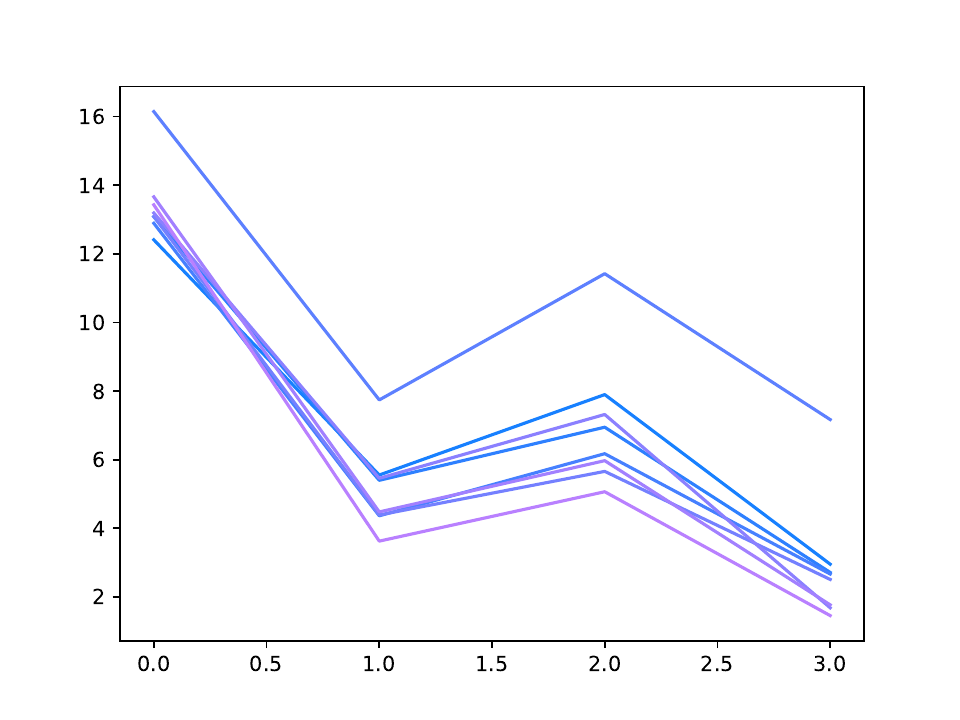}}
    \caption{Oversmoothing. Maximum node embedding distance after four main GNN operations. The plot is generated for MIS problem over random graphs with $N=200$. Oversmoothing is happening for graphs in (b) where the GNN training is unsuccessful.}
    \label{fig:smooth}
\end{figure}

A solution proposed in the literature for oversmoothing is graph sparsification by dropping edges \cite{rong2019dropedge}. This also explains how GNNs fail to learn good solutions on dense graphs as observed in Supplementary Figure \ref{fig:phase:r}. Motivated by \cite{rong2019dropedge}, we propose graph sparsification to improve the GNN's performance in \sys{}.

\subsection{Graph Sparsification}
With the observation that the inability of the GNNs to learn a good solution on some graphs is attributed to their high density, we explored the possibility of improving GNN learning through the process of graph sparsification. 
We experimented with the same graphs of Supplementary Figure \ref{fig:phase:r} with $N=200$ nodes and sparsified the ones that the GNN failed to learn a good solution by removing each edge with a given probability $Ps$. A similar method was conducted in \cite{rong2019dropedge} to alleviate the oversmoothing issue in graph neural networks. 
We show the performance  of \sys{} on the sparsified graphs with different values of $Ps$ in Supplementary Figure \ref{fig:phase:r:sparse}.
\begin{figure}[h!]
    \centering
 \includegraphics[width=0.45\textwidth]{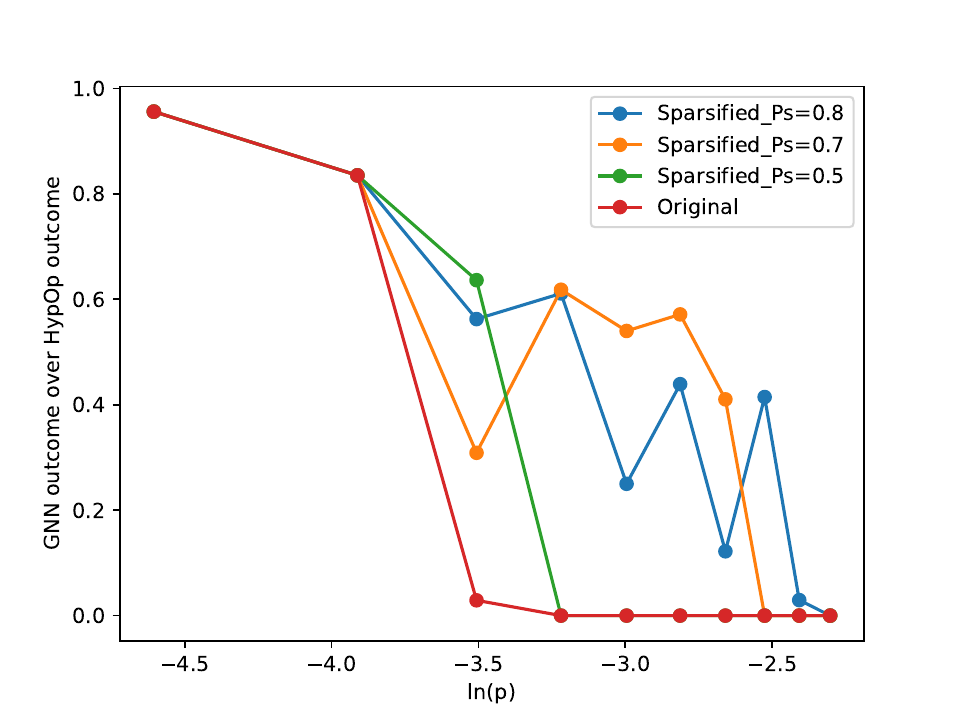}
    \caption{Graph Sparsification. The performance of PI-GNN over original and sparsified random graphs with $N=200$ for MIS problem.}
    \label{fig:phase:r:sparse}
\end{figure}
We can see in the figure that sparsification improves the performance of PI-GNN. However, there is still a gap between the overall outcome quality of \sys{} and the outcome of the  PI-GNN. 

Based on the findings presented in this section, we conclude that integrating Graph Neural Networks (GNNs) and HyperGNNs with other solvers yields significant advantages compared to employing them as standalone end-to-end solvers. This observation serves as additional motivation for the inclusion of the fine-tuning step using Simulated Annealing (SA) in \sys{}. This supplementary fine-tuning step ensures a reliable solution, even in cases where the GNN fails to learn a good solution.

  \section{Effect of Self-Loops in the HyperGCN}
 It is common to add self-loops to the graphs before training a GNN model on them. We observed that in \sys{} (and similarly in PI-GNN), adding self loops will deteriorate the performance of the GNN or HyperGNN.  The reason is most likely related to the results shown in the previous 
  section  where we showed that the GNNs fail to learn on graphs  that are denser than a given threshold. Adding self loops will make the graphs  denser than is tolerated by the GNN, and therefore, we did not add self loops and also set the diagonal  terms in the hypergraph adjacency matrix to 0 to create sparser graphs  and hypergraphs.

\break 

\end{document}